\documentclass[12pt]{article}
\usepackage{amsmath,amsfonts,amssymb,amsthm}

\newcommand{\al}{\alpha}
\newcommand{\del}{\delta}					

\newcommand{\lm}{\lambda}
\renewcommand{\th}{\theta}
\newcommand{\w}{\omega}
\newcommand{\1}{{\bf 1}}
\newcommand{\id}{\mbox{\rm id}}
\newcommand{\Hom}{\mbox{\rm Hom\,}}
\newcommand{\End}{\mbox{\rm End\,}}
\newcommand{\wt}{\mbox{\rm wt}}
\newcommand{\Res}{\mbox{\rm Res\,}}
\def\Z{\mathbb{Z}}
\def\Q{\mathbb{Q}}
\def\C{\mathbb{C}}
\def\N{\mathbb{N}}
\def\h{\mathfrak{h}}

\newcommand{\ots}{\otimes}
\newcommand{\ops}{\oplus}
\newcommand{\der}{\partial}
\newcommand{\Ops}{\bigoplus}

\def\cp{\cdotp}

\newcommand{\eqa}{\begin{eqnarray}}
\newcommand{\eeqa}{\end{eqnarray}}
\newcommand{\eqn}{\begin{eqnarray*}}
\newcommand{\eeqn}{\end{eqnarray*}}
\newcommand{\nn}{\nonumber}

\newtheorem{dfn}{Definition}[section]
\newtheorem{pro}[dfn]{Proposition}
\newtheorem{thm}[dfn]{Theorem}

\newtheorem{lem}[dfn]{Lemma}
\newtheorem{cor}[dfn]{Corollary}
\newtheorem{rem}[dfn]{Remark}

\makeatletter
\@addtoreset{equation}{section}

\makeatother

\def\bl{\begin{lem}\label}
\def\el{\end{lem}}
\def\bt{\begin{thm}\label}
\def\et{\end{thm}}
\def\bp{\begin{pro}\label}
\def\ep{\end{pro}}
\def\br{\begin{rem}\label}
\def\er{\end{rem}}
\def\bc{\begin{cor}\label}
\def\ec{\end{cor}}
\def\bd{\begin{dfn}\label}
\def\ed{\end{dfn}}
\def\proof{{\it Proof. }}
\renewcommand{\ss}[1]{\mbox{\scriptsize{${#1}$}}}
\def\B{\langle}
\def\K{\rangle}
\def\qed{{\mbox{$\ \ \ \square$}}}
\def\swt{\mbox{$\ss{\wt}$}}
\def\vs{\vspace{3mm}}
\def\hh{\hat{\mathfrak{h}}}

\def\Pa{\mbox{$\mathcal{P}_{L}$}}
\def\Qa{\mbox{$\mathcal{Q}_{L}$}}
\newcommand{\NO}{\,{\raise0.25em\hbox
{$\mathop{\hphantom {\cdot}}
\limits^{_{\circ}}_{^{\circ}}$}}\,}
\newcommand{\M}[1]{M(1)^{#1}}
\newcommand{\Ml}[1]{M(1,#1)}
\newcommand{\V}[1]{V_{L}^{#1}}
\newcommand{\Vl}[1]{V_{#1 +L}}
\newcommand{\Mt}[1]{M(1)(\th)^{#1}}
\newcommand{\Va}[1]{V_{\al/2+L}^{#1}}

\newcommand{\fusion}[3]{\mbox{${\ss{\mbox{$\left(\begin{array}
{cr}#3\\ #1\ \ \ #2 \end{array}\right)$}}}$}}

\begin{document}
\title{\Large  Fusion Rules for the Charge Conjugation
Orbifold}
\author{\large Toshiyuki Abe}
\date{\small\it Department of Mathematics, Graduate School of
Science, Osaka University,\\
 Toyonaka, Osaka 560-0043, Japan\\
{\rm e-mail: sm3002at@ecs.cmc.osaka-u.ac.jp}}
\maketitle
\begin{abstract}
We completely determine fusion rules for irreducible modules of the
charge conjugation orbifold.
\end{abstract}

\baselineskip 7mm

\section{Introduction}
The charge conjugation orbifold $\V{+}$ is the orbifold of the lattice vertex operator algebra $\V{}$ associated to a rank one even lattice $L$ by the automorphism $\th$
given by extending the $-1$-isometry of $L$ (cf. \cite[Section 6.1]{KT}). 
The set of all equivalence class of irreducible $\V{+}$-modules consists of $k+3$
modules derived from irreducible (untwisted) $\V{}$-modules (we call them
untwisted type modules) and $4$ modules from irreducible $\th$-twisted modules
(we call them twisted type modules) (see \cite{DN2}), where $k$ is the half square
length of the generator of $L$. In this paper we completely determine the fusion
rule for the irreducible $\V{+}$-modules.  The intertwining operators for
$\V{}$-modules constructed in \cite{DL} give rise to intertwining operators for
untwisted type modules.  We construct intertwining operators involving twisted
type modules by means of the twisted intertwining operators constructed in
\cite{FLM}.  The fusion rules and explicit forms of intertwining operators for the
free bosonic orbifold vertex operator algebra $\M{+}$ determined in \cite{A} play
important roles in analyzing intertwining operators for $\V{+}$.

The vertex operator algebra $\V{+}$ and its irreducible modules are constructed as
follows:  Let $L=\Z\al$ be a rank one even lattice with a $\Z$-bilinear form
$\B\cdotp,\cdotp\K$ defined by $\B \al,\al\K=2k$ for a positive integer $k$. Set
$\h=\C\ots_{\Z} L$ and extend the $\Z$-bilinear form to a $\C$-bilinear form on
$\h$ in the canonical way.  Let $\hh=\h\ots\C[t,t^{-1}]\ops\C K$ be its affinization
with the center $K$.  Then the Fock space $\M{ }=S(\h\ots t^{-1}\C[t^{-1}])$ is a
simple vertex operator algebra with central charge $1$.  Let
$\C[\h]=\ops_{\lm\in\h}\C e_{\lm}$ be the group algebra of the abelian group
$\h$, and set 
$\C[M]=\ops_{\lm\in M}\C e_{\lm}$ for a subset $M$ of $\h$. 
It is known that $\V{}=\M{}\ots\C[L]$ is a simple vertex operator algebra with
central charge $1$ (cf. \cite{FLM}), and $\Vl{\lm}=\M{}\ots\C[\lm+L]$ is an
irreducible $\V{}$-module for all $\lm\in L^{\circ}$, where $L^{\circ}$ is the dual
lattice of $L$.  Moreover all irreducible $\V{}$-modules are given by the set
$\{\,\Vl{\lm}\,|\,
\lm+L\in L^{\circ} /L\,\}$ (cf. \cite{D1}).
Let $\th$ be the involution of $L$ defined by $\th(\beta)=-\beta$ for $\beta\in L$.
Then the involution $\th$ can be lifted to an isomorphism of $V_{L^{\circ}}$, and
the $\th$-invariant subspace of $\V{}$ becomes a simple vertex operator algebra
with central charge $1$, denoted by $\V{+}$. The automorphism $\th$ induces an
$\V{+}$-module isomorphism  from $\Vl{\beta}$ to $\Vl{-\beta}$ for $\beta \in
L^{\circ}$.  For a $\th$-invariant subspace $W$ of $V_{L^{circ}}$, we denote the
$\pm1$-eigenspaces by $W^{\pm}$ respectively.  Then $\V{\pm}$, $\Va{\pm}$
and $\Vl{r\al/2k}$ for $1\leq r\leq k-1$ are irreducible $\V{+}$-modules (see
\cite{DN2}). 

Let $\hh[-1]=\h\ots t^{1/2}\C[t,t^{-1}]\ops\C K$ be the twisted
affine Lie algebra and set $\Mt{}=S(\h\ots t^{-1/2}\C[t^{-1}])$. Then
$\Mt{}$ is a unique irreducible $\th$-twisted $\M{}$-module (cf. \cite{FLM} and
\cite{D2}). The automorphism $\th$ acts on $\Mt{}$, and the $\pm1$-eigenspaces
$\Mt{\pm}$ become irreducible $\M{+}$-modules (see \cite{DN1}).  Let $T^{1}$ and
$T^{2}$ be irreducible $\C[L]$-modules on which $e_{\al}$ acts $1$  and $-1$
respectively. Then the tensor products $\V{T^{i}}=\Mt{}\ots T^{i}$ ( $i=1,2$) are
irreducible $\th$-twisted $\V{}$-modules, and their
$\pm1$-eigenspaces $\V{T^{i},\pm}$ for $\th$ become irreducible $\V{+}$-modules
(\cite{DN2}).  In \cite{DN2}, it is proved that every 
irreducible $\V{+}$-module is isomorphic to one of the irreducible modules
$\V{\pm},\Va{\pm},\Vl{r\al/2k}$ for $1\leq r\leq k-1$ and $\V{T_{i},\pm}$ for
$i=1,2$.

For a vertex operator algebra $V$ and its modules $W^{1},W^{2}$ and $W^{3}$,
 the dimension of the vector space $I_{V}\fusion{W^{1}}{W^{2}}{W^{3}}$ of all
intertwining operators of type $\fusion{M^{1}}{M^{2}}{M^{3}}$ is called the fusion
rule of corresponding type and denoted by $N_{W^{1}W{2}}^{W^{3}}$. It is known
that fusion rules have the following symmetry;
\eqa I_{V}\fusion{W^{1}}{W^{2}}{W^{3}}\cong\fusion{W^{2}}{W^{1}}{W^{3}}
\cong\fusion{W^{1}}{(W^{3})'}{(W^{2})},\label{SY1}\eeqa
where $W'$ means the contragredient module of $W$ (see \cite{FHL,HL}). 

We give the correspondence of irreducible $\V{+}$-modules and contragredient
modules (see Proposition \ref{Con1}).  The correspondence and the symmetry of
fusion rules (\ref{SY1}) are very useful in reducing the arguments to determine the
fusion rules for $\V{+}$. 

We explain the method of determining the fusion rules for $\V{+}$ in more detail. 
Let $W^{1}, W^{2}$ and $W^{3}$ be $\V{+}$-modules, and suppose that $W^{1}$ and
$W^{2}$ contain $\M{+}$-submodules $M^{1}$ and $M^{2}$ respectively.  
Then we have a canonical restriction map 
$$I_{\V{+}}\fusion{W^{1}}{W^{2}}{W^{3}}\to I_{\M{+}}\fusion{M^{1}}{M^{2}}{W^{3}},\ 
\mathcal{Y}\mapsto\mathcal{Y}|_{N^{1}\ots N^{2}}.$$
It is known that if $W^{1}$ and $W^{2}$ are irreducible, the restriction map is
injective (cf. \cite{DL}, Proposition 11.9). Therefore we then  have
\eqa
\dim I_{\V{+}}\fusion{W^{1}}{W^{2}}{W^{3}}\leq
\dim I_{\M{+}}\fusion{M^{1}}{M^{2}}{W^{3}}.\label{IJ1}
\eeqa
We also prove that all irreducible $\V{+}$-modules are
completely reducible as $\M{+}$-modules and that the multiplicity of each
irreducible $\M{+}$-module is at most one.  
Using this fact, (\ref{IJ1}) and fusion
rules for $\M{+}$, we show that fusion rules for $\V{+}$ are zero or one. 
The formula (\ref{IJ1}) also shows that for irreducible
$\V{+}$-modules $W^{1},W^{2}$ and $W^{3}$,  if there are $\M{+}$-submodules
$M^{1}$ of $W^{1}$ and $M^{2}$ of $W^{2}$ such that the fusion rule
$N_{M^{1}M^{2}}^{W^{3}}$ for $\M{+}$ is zero,  then the fusion rule
$N_{W^{1}W^{2}}^{W^{3}}$ for $\V{+}$ is zero. 
For almost of irreducible $\V{+}$-modules
$W^{1},W^{2}$ and $W^{3}$ for which the fusion rule
$N_{W^{1}W^{2}}^{W^{3}}$ is zero, we can find such $\M{+}$-submodules $M^{1}$ of
$W^{1}$ and $M^{2}$ of $W^{2}$.  But there are irreducible modules $W^{1},W^{2}$
and $W^{3}$ such that the fusion rule $N_{W^{1}W^{2}}^{W^{3}}$ is zero although
the fusion rule $N_{M^{1}M^{2}}^{W^{3}}$ is nonzero for any $\M{+}$-submodules
$M^{1}$ of $W^{1}$ and $M^{2}$ of $W^{2}$ (for example
$W^{1}=\V{-}$, $W^{2}=W^{3}=\Va{+}$). 
In such case, to show that intertwining operators $\mathcal{Y}$ of corresponding
type are zero, we choose irreducible $\M{+}$-submodules $M^{1}$ and $M^{2}$
from $W^{1}$ and $W^{2}$ respectively.  
Since the fusion rule $N_{M^{1}M^{2}}^{W^{3}}$ is nonzero
and $W^{3}$ is a direct sum of irreducible $\M{+}$-modules as $W^{3}=\ops_{i\in
I}M^{3}_{i}$,  we see that the restriction of $\mathcal{Y}$ to $M^{1}\ots M^{2}$ is
a linear combination of intertwining operators of types
$\fusion{M^{1}}{M^{2}}{M^{3}_{i}}$.  Then the explicit forms of intertwining
operators of types $\fusion{M^{1}}{M^{2}}{M^{3}_{i}}$ shows $\mathcal{Y}=0$.
 
The nonzero fusion rules is provided by constructing nontrivial intertwining
operators explicitly. The constructions is separated in two cases; one is the case all
modules are of untwisted types, and other is the case some modules are of twisted
types. 

In the case all modules are of untwisted types, the nontrivial intertwining
operators are essentially given in \cite{DL}. 
\cite{DL} construct a nontrivial intertwining operator $\mathcal{Y}_{\lm\mu}$ for
$\V{}$ of type
$\fusion{\Vl{\lm}}{\Vl{\mu}}{\Vl{\lm+\mu}}$ for $\lm,\mu\in
L^{\circ}$. 
The intertwining operator $\mathcal{Y}_{\lm\mu}$ gives rise to a nonzero
intertwining operator for $\V{+}$ of type
$\fusion{\Vl{\lm}}{\Vl{\mu}}{\V{\lm+\mu}}$. 
Since $\th$ induces a $\V{+}$-module isomorphism from $\Vl{\lm}$ to
$\Vl{-\lm}$ for $\lm\in L^{\circ}$, the operator
$\mathcal{Y}_{\lm,-\mu}\circ{\th}$ defined by
$\mathcal{Y}_{\lm,-\mu}\circ{\th}(u,z)v=\mathcal{Y}_{\lm,-\mu}(u,z)\th(v)$ for
$u\in\Vl{\lm}$ and $v\in\Vl{\mu}$ gives a nonzero intertwining operator of type
$\fusion{\Vl{\lm}}{\Vl{\mu}}{\V{\lm-\mu}}$. 
Then all nonzero intertwining
for untwisted type modules are given by restricting $\mathcal{Y}_{\lm\mu}$ or
$\mathcal{Y}_{\lm,-\mu}\circ\th$ to irreducible $\V{+}$-modules.  

In the case some modules are of twisted types, we construct nonzero intertwining
operators  as follows:  In \cite{A}, an intertwining operator
$\mathcal{Y}^{\th}$ for $\M{+}$ of type $\fusion{\Ml{\lm}}{\Mt{}}{\Mt{}}$ for
$\lm\in L^{\circ}$ is constructed following \cite{FLM}. 
As in \cite{DL}, for $\lm\in L^{\circ}$, we give an linear isomorphism $\psi_{\lm}$
of $ T^{1}\ops T^{2}$ which satisfies
$e_{\al}\psi_{\lm}=(-1)^{\B\al,\lm\K}\psi_{\lm}e_{\al}=\psi_{\lm+\al}$, and define
$\tilde{\mathcal{Y}}$ by
$\tilde{\mathcal{Y}}(u,z)=\mathcal{Y}^{\th}(u,z)\ots\psi_{\gamma}$ for
$\gamma\in\lm+L$ and
$u\in\Ml{\gamma}$. Then for $\lm\in L^{\circ}$ and $i,j=1,2$ which satisfy
$(-1)^{\B\lm,\al\K+\del_{i,j}+1}=1$, $\tilde{\mathcal{Y}}$ gives rise to an
intertwining operator of type $\fusion{\Vl{\lm}}{\V{T_{i}}}{\V{T_{j}}}$, and all
nonzero intertwining operators in this case are given by restricting
$\tilde{\mathcal{Y}}$ to irreducible $\V{+}$-modules and by
using symmetry of fusion rules (\ref{SY1}).  

The organization of this paper is as follows: 
We recall definitions of modules for a vertex operator algebra and fusion rules in
Section \ref{S2.1},   we review the vertex operator algebras $\M{+}$ and $\V{+}$
and their irreducible modules in Section \ref{S3.1}.  In Section \ref{S2.4} we state
the fusion rules for $\M{+}$, and discuss the contragredient modules for $\V{+}$.  
In Section \ref{S3.2}, we give the  irreducible decompositions of irreducible
$\V{+}$-modules as $\M{+}$-modules and prove that the fusion rules for $\V{+}$
are zero or one.  In Section \ref{S4.4}, we state the main theorem (Theorem
\ref{T4}).  In Section \ref{S4.1} and \ref{S4.2}, the proof of the
main theorem is given. 
In Section \ref{S4.1}, we determine the fusion rule
$N_{W^{1}W^{2}}^{W^{3}}$ for untwisted type modules $W^{i}$ $(i=1,2,3)$.
In Section \ref{S4.2}, fusion rule of type $\fusion{W^{1}}{W^{2}}{W^{3}}$
are determine in the case some of $W^{i}$ $(i=1,2,3)$ are twisted type module.


\section{Preliminaries}\label{S2}
In Section \ref{S2.1}, we recall the definition of a $g$-twisted module for a 
vertex operator algebra and its automorphism $g$ of finite order and that of an
intertwining operator following \cite{FLM,FHL,DMZ} and \cite{DLM1}. In Section
\ref{S3.1}, we review constructions of vertex operator algebras $\M{+}$, $\V{+}$ and their
irreducible modules following \cite{FLM,DL,DN1,DN2}. In Section \ref{S2.4}, we state the fusion
rules for $\M{+}$ obtained in \cite{A} (see Theorem \ref{T3}) and discuss the contragredient modules for $\V{+}$. 

Throughout this paper, $\N$ is the set of
nonnegative integers and $\Z_{+}$ is the set of positive integers.

\subsection{Modules, intertwining operators and fusion rules}\label{S2.1}
Let $(V,Y,\1,\w)$ be a vertex operator algebra and $g$ an automorphism of $V$ of
order $T$. Then $V$ is decomposed into the direct sum of eigenspaces for $g$:
\eqa
V=\Ops_{r=0}^{T-1} V^{r},\ V^{r}=\{ \ a\in V\ |\ g(a)=e^{-\frac{2\pi i
r}{ T}} a\ \}\nn.
\eeqa
 
A {\em $g$-twisted
$V$-module} is a $\C$-graded vector space $M=\ops_{\,\lm\in\C\,} M(\lm)$ such that each $M(\lm)$ is finite dimensional and for fixed
$\lm\in\C$, $M(\lm+{n/T})=0$ for sufficiently small
integer $n$, and equipped with a linear map
\eqa
Y_{M}:V&\to&(\End M)\{z\},\nn\\
a&\mapsto&Y_{M}(a,z)=\sum_{n\in\Q} a_{n}^{M}z^{-n-1},\
(a_{n}^{M}\in \End M)\nn
\eeqa
such that the following conditions hold for $0\leq r\leq T-1,\ a\in
V^{r},b\in V\ and\ u\in M$:
\eqn
Y_{M}(a,z)=\sum_{n\in {r/ T}+\ss{\Z}} a_{n}^{M} z^{-n-1},\
Y_{M}(a,z)v\in z^{-{\frac{r}{T}}}M((z)),
\eeqn
\eqa
&&z_{0}^{-1}\delta\left(
{\frac{z_{1}-z_{2}}{z_{0}}}\right)Y_{M}(a,z_{1})Y_{M}(b,z_{2})-z_{0}^{-1}\delta
\left({\frac{z_{2}-z_{1}}{ -z_{0}}}\right)Y_{M}(b,z_{2})Y_{M}(a,z_{1})\nn\\
&&{}=z_{2}^{-1}\delta\left(
{\frac{z_{1}-z_{0}}{z_{2}}}\right)\left(
{\frac{z_{1}-z_{0}}{z_{2}}}\right)^{-{\frac{r}{T}}}Y_{M}(Y(a,z_{0})b,z_{2}),\nn
\eeqa
$$ Y_{M}(\1,z)=\id_{M},$$
\eqn 
L(0)v=\lm v\hbox{ for $v\in M(\lm)$, }
\eeqn
where we set $Y_{M}(\w,z)=\sum_{n\in\Z}L(n)z^{-n-2}$.

A $g$-twisted $V$-module is denoted by $(M,Y_{M})$ or
simply by $M$. In the case $g$ is the identity of $V$, a $g$-twisted
$V$-module is called a {\it $V$-module}. An element $u\in M(\lm)$ is called a {\it homogeneous element of weight}
$\lm$. We denote the weight by $\lm = \wt(u)$. We write the
component operator $a_{n}^{M}\ (a\in V,\ n\in\Q)$ by $a_{n}$ for simplicity.

For a $V$-module $M$, it is known that the restricted
dual $M'=\ops_{\lm\in\C} M(\lm)^{*}$ with the vertex operator $Y_{M}^{*}(a,z)$ for $a\in V$ defined by 
\eqa\B Y_{M}^{*}(a,z)u',v\K =\B
u',Y_{M}(e^{z L(1)} (-z^{-2})^{L(0)}a,z^{-1})v\K\nn\eeqa 
for $u'\in M',\ v\in M$ is a $V$-module (cf. \cite{FHL}).
The $V$-module $(M',Y_{M}^{*})$ is called the {\it contragredient
module of} $M$. The double contragredient module $(M')'$ of $M$ is
naturally isomorphic to $M$, and therefore if $M$ is irreducible, then $M'$ is also irreducible (see \cite{FHL}). 

\bd{D3.7} Let $V$ be a vertex operator algebra and $(M^{i}, Y_{M^{i}})$
$(i=1,2,3)$ be $V$-modules. An {\rm intertwining operator for} $V$ {\rm of
type}
\fusion{M^{1}}{M^{2}}{M^{3}}
 is a linear map $\mathcal{Y}:M^{1}\ots M^{2}\to M^{3}\{z\}$, or equivalently,
\eqa
\mathcal{Y}:M^{1}&\to&(\Hom (M^{2},M^{3}))\{z\},\nn\\
v&\mapsto&\mathcal{Y}(v,z)=\sum_{n\in\C}v_{n} z^{n}\
(v_{n}\in\Hom(M^{2},M^{3}))\nn
\eeqa
such that for $a\in V, v\in M^{1}\ and\ u\in M^{2}$,
following conditions are satisfied:

For fixed $n\in\C$, $v_{n+k}u=0$ for sufficiently large integer $k$,

\eqa
&&z_{0}^{-1}\delta\left(
{\frac{z_{1}-z_{2}}{z_{0}}}\right)Y_{M^{3}}(a,z_{1})\mathcal{Y}
(v,z_{2})-z_{0}^{-1}\delta\left({\frac{z_{2}-z_{1}}{
-z_{0}}}\right)\mathcal{Y}(v,z_{2})Y_{M^{2}}(a,z_{1})\nn\\
&&{}=z_{2}^{-1}\delta\left( {\frac{z_{1}-z_{0}}{z_{2}}}\right)
\mathcal{Y}(Y_{M^{1}}(a,z_{0})v,z_{2}),\label{JI1}
\eeqa
\eqn
{\frac{d}{dz}} \mathcal{Y}(v,z)= \mathcal{Y}(L(-1)v,z).\label{DP1}\eeqn
\ed
The vector space of all intertwining operators of type \fusion{M^{1}}{M^{2}}{M^{3}} is denoted by
$I_{V}$\fusion{M^{1}}{M^{2}}{M^{3}}. The dimension of the vector
space $I_{V}\fusion{M^{1}}{M^{2}}{M^{3}}$ is called the {\it
fusion rule} of corresponding type and denoted by $N_{M^{1}M^{2}}^{M^{3}}$.
Fusion rules have the following symmetry (see
\cite{FHL} and \cite{HL}).
\bp{P9} Let $M^{i}$ $(i=1,2,3)$ be $V$-modules. Then there exist natural
isomorphisms 
\eqn
I_{V}\fusion{M^{1}}{M^{2}}{M^{3}}\cong I_{V}\fusion{M^{2}}{M^{1}}
{M^{3}}\hbox{ and }I_{V}\fusion{M^{1}}{M^{2}}{M^{3}}\cong
I_{V}\fusion{M^{1}}{(M^{3})'}{(M^{2})'}.
\eeqn
\ep
The following lemma is often used in later sections.
\bl{P1}{\rm(\cite{DL})} Let $V$ be a vertex operator algebra, and let $M^{1}$  and
$M^{2}$ be irreducible $V$-modules and $M^{3}$ a $V$-module. If $\mathcal{Y}$
is a nonzero intertwining operator of type $\fusion{M^{1}}{M^{2}}{M^{3}}$, then
$\mathcal{Y}(u,z)v\neq0$  for any nonzero vectors $u\in M^{1}$ and $v\in M^{2}$.
\el
As a direct consequence of Lemma \ref{P1}, we have
\bc{C8} Let $V, M^{i}$ $(i=1,2,3)$ be as in Lemma \ref{P1}, and let $U$ be a
vertex operator subalgebra of $V$ with same Virasoro element, $N^{i}$ a
$U$-submodule of
$M^{i}$ for $i=1,2$. Then the restriction map 
\eqn
 I_{V}\fusion{M^{1}}{M^{2}}{M^{3}}\to
I_{U}\fusion{N^{1}}{N^{2}}{M^{3}},\ \mathcal{Y}\mapsto\mathcal{Y}|_{N^{1}\ots N^{2}},
\eeqn 
is injective. In particular, we have
\eqa\dim I_{V}\fusion{M^{1}}{M^{2}}{M^{3}}\leq
\dim I_{U}\fusion{N^{1}}{N^{2}}{M^{3}}.\nn
\eeqa
\ec
Let $V, M^{i}$ $(i=1,2,3)$, $U$ and  $N^{i}$ $(i=1,2)$ be as in Corollary \ref{C8}.
Suppose that $M^{3}$ is decomposed into a direct sum of irreducible $U$-modules
as $M^{3}=\ops_{i} L^{i}$. Then there is  an isomorphism
\eqn I_{U}\fusion{N^{1}}{N^{2}}{\ops_{i} L^{i}}\cong \ops_{i}
I_{U}\fusion{N^{1}}{N^{2}}{L^{i}}.
\eeqn
Therefore by Corollary \ref{C8}, we have an in equality
\eqa
\dim I_{V}\fusion{M^{1}}{M^{2}}{M^{3}}\leq \sum_{i}\dim
I_{U}\fusion{N^{1}}{N^{2}}{L^{i}}.\label{Ineq2}
\eeqa

Another consequence of Lemma \ref{P1} is
\bl{Le1}Let $V$ be a simple vertex operator algebra, and let $M^{1}$ and $M^{2}$
be  irreducible $V$-modules. If the fusion rule of type $\fusion{V}{M^{1}}{M^{2}}$
is nonzero, then $M^{1}$ and $M^{2}$ are isomorphic to each other as
$V$-modules.\el
\proof Let $\mathcal{Y}$ be an intertwining operator of type
$\fusion{V}{M^{1}}{M^{2}}$.  Consider the operator ${\mathcal  Y}(\1,z)$. By the
$L(-1)$-derivative property (\ref{DP1}), we see that ${\mathcal Y}(\1,z)$ is
independent on $z$. Denote
$f={\mathcal Y}(\1,z)\in\Hom(M^{1},M^{2})$. Since $V$ is simple and $M^{1}$ is
irreducible, Proposition \ref{P1} implies that $f$ is nonzero. By Jacobi identity
(\ref{JI1}), we have a commutation relation 
\eqn 
[a_{n},\mathcal{Y}(\1,z)]=\sum_{i=0}^{\infty}{\mbox{$\left(\begin{array}{c}n\\i
\end{array}\right)$}}\mathcal{Y}(a_{i}\1,z) z^{n-i}=0
\eeqn
for $a\in V$ and $n\in\Z$. Hence $f$ is a nonzero $V$-module homomorphism from
$M^{1}$ to $M^{2}$. Since $M^{1}$ and $M^{2}$ are irreducible, $f$ is in fact
isomorphism. Therefore $M^{1}$ is isomorphic to $M^{2}$.\qed


\subsection{Vertex operator algebra $\V{+}$ and its irreducible
modules}\label{S3.1}
We discuss the constructions of vertex operator algebras $\M{}$, $\V{}$ and their 
irreducible (twisted) modules following \cite{FLM, DL,D1} and \cite{D2}. We also
refer to the vertex operator algebras $\M{+}$, $\V{+}$ and irreducible modules for
them classified in \cite{DN1,DN2}.
  
Let $L$ be an even lattice of rank
$1$ with a nondegenerate positive definite ${\Z}$-bilinear
form $\B\cdotp,\cdotp\K$, and $\h=\C\ots_{\Z} L$. Then $\h$ has the
nondegenerate symmetric $\C$-bilinear form given by extending the form
$\B\cp,\cp\K$ of $L$. Let $\C[\h]$ be the group algebra of $\h$ with a basis
$\{\,e_{\lm}\,|\,\lm\in\h\,\}$. For a subset $M$ of
$\h$, set $\C[M]=\ops_{\lm\in M}\C e_{\lm}$. 

Let $\hh =\h\ots \C[t,t^{-1}]\ops\C K$ be a Lie algebra with
the commutation relation given by $[X\ots t^{m},X'\ots
t^{n}]=m\,\del_{m+n,0}\,\B X,X'\K\, K,[K,\hh]=0$ for $X,X'\in\h$ and
$m,n\in\Z$. Then $\hh^{+}=\h\ots
\C[t]\ops \C K$ is a subalgebra of $\hh$, and the group algebra
$\C[\h]$ becomes a $\hh^{+}$-module by the action  $\rho(
X\ots t^{n}) e_{\lm}=\del_{n,0}\B X, \lm\K e_{\lm}$ and $\rho(K)
e_{\lm}=e_{\lm}$ for
$\lm\in \h$, $X\in \h$ and $n\in\N$. It is clear that for a subset $M$ of
$\h$ the subspace $\C[M]$ is a $\hh^{+}$-submodule of $\C[\h]$. Set $V_{M}$ the induced
module of $\hh$ by $\C[M]$:
\eqn 
V_{M}=U(\hh)\ots_{U(\hh^{+})} \C[M]\cong S(\h\ots
t^{-1}\C[t^{-1}])\ots\C[M]\hbox{ ({\rm linearly}),}
\eeqn 
where $U(\cp)$ means a universal enveloping algebra. Denote the
action of
$X\ots t^{n}$ ($X\in\h, n\in\Z$) on $V_{\h}$  by $X(n)$ and set
$X(z)=\sum_{n\in\Z}X(n) z^{-n-1}$ for $X\in\h$. For
$\lm\in \h$, the vertex operator associated with $e_{\lm}$ is
defined by
\eqa
\mathcal{Y}^{\circ}(e_{\lm},z)=\exp\left(\sum_{n=1}^{\infty}
{\frac{\lm(-n)}{n}}z^{n}\right)
\exp\left(-\sum_{n=1}^{\infty}{\frac{\lm(n)}{n}}z^{-n}\right)
e_{\lm}z^{\lm(0)},\label{O1}
\eeqa 
where $e_{\lm}$ in the right-hand side means the left
multiplication of $e_{\lm}\in \C[\h]$ on the group algebra $\C[\h]$,  and $z^{\lm(0)}$ is an
operator on $V_{\h}$ defined by $z^{\lm(0)}u=z^{\B \lm,\mu\K}u$ for $\mu\in \h$ and 
$u\in\Ml{\mu}$. For $v=X_{1}(-n_{1})\cdots
X_{m}(-n_{m})\, e_{\lm}\in V_{\h}$ ($X_{i}\in\h$ and
$n_{i}\in\Z_{+}$), the corresponding vertex operator is defined by
\eqa \mathcal{Y}^{\circ}(v,z)=\NO\der^{(n_{1}-1)}X_{1}(z)\cdots\der^{(n_{m}-1)}
X_{m}(z) \mathcal{Y}^{\circ}(e_{\lm},z)\NO,\label{O2}
\eeqa
where $\der^{(n)}=(\frac{1}{n!})(d/dz)^{n}$, and the normal ordering
$\NO\cp\NO$  is an operation which reorders so that $X(n)$ ($X\in\h,n< 0$) and 
$e_{\lm}$ to be placed to the left of $X(n)$ ($X\in\h,n\geq 0$) and $z^{\lm(0)}$.
We extend $\mathcal{Y}^{\circ}$ to $V_{\h}$ by linearity. We denote 
$Y(a,z)=\mathcal{Y}^{\circ}(a,z)$ when $a$ is in $\V{}$.

Set $L=\Z\al$ and $\B\al,\al\K=2 k$ for $k\in\Z_{+}$, and $L^{\circ}=\{\,\lm\in\h\,|\,\B\lm,\al\K\in\Z\,\}$, the dual lattice of $L$. 
Let $h=\al /\sqrt{2k}$ be the orthonormal basis of $\h$ and set $\1=1\ots e_{0}$ 
and $\w=(1/2)\,h(-1)^{2}e_{0}$. Then $(\V{},Y,\1,\w)$
is a simple vertex operator algebra with central charge $1$ and  for $\lm+L\in
L^{\circ}/L$, $(V_{\lm+L},Y)$ is an irreducible module for $\V{}$.
Furthermore $V_{\lm+L}$ for $\lm+L\in L^{\circ}/L$ give all
inequivalent irreducible $\V{}$-modules. Set $\M{}=S(\h\ots
t^{-1}\C[t^{-1}])\ots e_{0}\subset
\V{}$, then ($\M{},Y,\1,\w)$ is a simple vertex operator algebra. If we set $\Ml{\lm}=U(\hh)\ots_{U(\hh^{+})}\C e_{\lm}$ for
each $\lm\in\h$, then $(\Ml{\lm},\mathcal{Y}^{\circ})$ becomes an irreducible 
$\M{}$-module (see \cite{D1,DL}).

Let $\th$ be a linear isomorphism of $V_{\h}$ defined by 
\eqa 
\th(X_{1}(-n_{1})X_{2}(-n_{2})\cdots X_{\ell}(-n_{\ell})\ots
e_{\lm})=(-1)^{\ell}\,X_{1}(-n_{1})X_{2}(-n_{2})\cdots 
X_{\ell}(-n_{\ell})\ots e_{-\lm},\nn
\eeqa
for $X_{i}\in\h,n\in\Z_{+}$ and $\lm\in \h$. Then $\th$
induces automorphisms of $\V{}$ and $\M{}$. For a $\th$-invariant
subspace $W$ of $V_{\h}$, we denote the $\pm 1$-eigenspaces of
$W$ for $\th$ by $W^{\pm}$. Then $(\V{+},Y,\1,\w)$ and $(\M{+},Y,\1,\w)$ are 
vertex operator algebras. Furthermore $\M{\pm}$ and $\Ml{\lm}$ for $\lm\neq0$
are irreducible $\M{+}$-modules, and $\th$ induces an $\M{+}$-module
isomorphism between $\Ml{\lm}$and $\Ml{-\lm}$ (see\cite{DN1}). As to 
$\V{+}$-modules, $\V{\pm}$, $V_{{\al}/{2}+L}^{\pm}$ and $V_{r\al/2k+L}$ for
$1\leq r\leq k-1$ are irreducible modules (see \cite{DN2}) and $\th$ induces a
$\V{+}$-module isomorphism between $\Vl{\lm}$ and $\Vl{-\lm}$ for
$\lm\in L^{\circ}$.

Now we review the construction of $\th$-twisted $\V{}$-modules following
\cite{FLM,D2}. Let
$\hh[-1]=\h\ots t^{1/2}\C[t,t^{-1}]\ops \C K$ be a Lie algebra with
commutation relation $[X\ots t^{m},X'\ots t^{n}]=m\,\del_{m+n,0}\,\B
X,X'\K\, K,$ $[K,\hh[-1]]=0$ for $X,X'\in\h$ and $m,n\in 1/2+\Z$. Then $\C$
becomes a one-dimensional module for $\hh[-1]^{+}=\h\,\ots
t^{1/2}\,\C[t]\ops \C K$ by defining the actions by $\rho(
X\ots t^{n}) 1=0$ and $\rho(K)\,1=1$ for $X\in \h$ and
$n\in 1/2+\N$. Set $\Mt{}$ the induced $\hh[-1]$-module :
\eqn
 \Mt{}=U(\hh[-1])\ots_{U(\hh[-1]^{+})}\,\C\cong
S\left(\h\ots t^{-\frac{1}{2}}\C[t^{-1}]\right)\ \ ({\rm linearly}).
\eeqn
Denote the action of $X\ots t^{n}$ ($X\in\h, n\in 1/2+\Z$) on $\Mt{}$ by
$X(n)$, and set $X(z)=\sum_{n\in 1/2+\Z} X(n)\, z^{-n-1}$.  For $\lm\in
L^{\circ}$ a twisted vertex operator associated with
$e_{\lm}\in V_{\h}$ is defined by
\eqa \mathcal{Y}^{\th}(e_{\lm},z)=2^{-\B\lm,\lm\K}z^{-\frac{\B
\lm,\lm\K}{2}}\exp\left(\sum_{n\in1/2+\N} {\frac{\lm(-n)}{n}}z^{n}\right)
\exp\left(-\sum_{n\in1/2+\N}{\frac{\lm(n)}{n}}z^{-n}\right).\label{U1}
\eeqa
For $v=X_{1}(-n_{1})\cdots
X_{m}(-n_{m})e_{\lm}\in V_{L^{\circ}}$ ($X_{i}\in\h$ and
$n_{i}\in \Z_{+}$), set
\eqa W^{\th}(v,z)=\NO\der^{(n_{1}-1)}X_{1}(z)\cdots\der^{(n_{m}-1)}
X_{m}(z) \mathcal{Y}^{\th}(v_{\lm},z)\NO,\label{U2}
\eeqa
and extend it to $V_{L^{\circ}}$ by linearity, where the normal ordering
$\NO\cp\NO$  is an operation which reorders so that $X(n)$ ($X\in\h,n<0$) to be
placed to the left of $X(n)$ ($X\in\h,n> 0$). Let $c_{mn}\in\Q$ be coefficients
defined by the formal power series expansion
\eqn
\sum_{m,n\geq 0}c_{mn}x^{m}y^{n}=-\log
\left(\frac{(1+x)^{\frac{1}{2}}+(1+y)^{\frac{1}{2}}}{2}\right),
\eeqn
and set $\Delta _{z}=\sum_{m,n\geq 0} c_{mn}h(m)h(n) z^{-m-n}$. Then
the twisted vertex operator associated to $u\in V_{L^{\circ}}$ is defined by
\eqa
\mathcal{Y}^{\th}(u,z)=W^{\th}(\exp(\Delta_{z})u,z).\label{U3}
\eeqa
If we write $Y^{\th}(a,z)=\mathcal{Y}^{\th}(a,z)$ for $a\in \M{}$,
the pair $(\Mt{},Y^{\th})$ is the unique irreducible $\th$-twisted $\M{}$-module.
 
Let $T_{1}$ and $T_{2}$ be irreducible $\C[L]$-modules which $e_{\al}$
acts as $1$ and $-1$ respectively, and set $\V{T_{i}}=\Mt{}\ots_{\C} T_{i}$
for $i=1,2$. For $u\in\Ml{\beta}$ ($\beta\in L$), the corresponding twisted vertex operator is defined by 
$Y^{\th}(u,z)=\mathcal{Y}^{\th}(u,z)\ots
e_{\beta}$. We extend $Y^{\th}$ to $\V{}$ by linearity.  Then $(\V{T_{i}},Y^{\th})$
($i=1,2$) are irreducible $\th$-twisted
$\V{}$-modules. Note that $\V{T_{i}}$ has a $\th$-twisted $\M{}$-module
structure. Let $t_{i}$ be a basis of $T_{i}$ for
$i=1,2$. Then we have a canonical $\th$-twisted $\M{}$-module isomorphism 
\eqa 
\phi_{i}:\Mt{}\to\V{T_{i}}:u\mapsto u\ots t_{i}.\label{ISO2}
\eeqa

The action of the automorphism $\th$ on $\Mt{}$ is defined by
\eqn \th(X_{1}(-n_{1})\cdots X_{m}(-n_{m})1)=(-1)^{m}X_{1}(-n_{1})\cdots
X_{m}(-n_{m})1, \eeqn
for $X_{i}\in\h, n_{i}\in 1/2+\N$. Set $\Mt{\pm}$  the $\pm 1$-eigenspaces
of $\Mt{}$ for $\th$ and $\V{T_{i},\pm}$ the $\pm 1$-eigenspaces of $\V{T_{i}}$
for $\th\ots 1$. Then $\Mt{\pm}$ and $\V{T_{i},\pm}$ ($i=1,2$) become irreducible
$\M{+}$-modules and irreducible $\V{+}$-modules respectively (see
\cite{DN1} and \cite{DN2} resp.). All irreducible $\M{+}$-modules and all irreducible $\V{+}$-modules are
classified in \cite{DN1} and
\cite{DN2} respectively.

\bt{T1}{\rm (1)\ (\cite{DN1})} The set
\eqa\{\M{\pm},\Mt{\pm},\Ml{\lm}(\cong \Ml{-\lm})\,|\,\lm\in\h-\{0\}\,\}
\eeqa
gives
all inequivalent irreducible $\M{+}$-modules.

\noindent
{\rm (2)\ (\cite{DN2})} The set
\eqa\{\V{\pm},\Va{\pm},\V{T_{i},\pm},\Vl{r\al/2k}\,
|\,i=1,2,1\leq r\leq k-1\,\}\label{IM1}
\eeqa
gives all inequivalent irreducible
$\V{+}$-modules.
\et
We call irreducible modules $\V{\pm}$, $\Va{\pm}$ and
$\Vl{r\al/2k}$ {\it untwisted type} modules, and call $\V{T_{i},\pm}$ $(i=1,2)$ 
{\it twisted type} modules. Here and further we write
$\lm_{r}=r\al/2k$ for $r\in\Z$.

\subsection{Fusion rules for $\M{+}$ and contragredient modules for
$\V{+}$}\label{S2.4} First we list up the fusion rules for $\M{+}$ determined in
\cite{A}. The fusion rules play central roles in determining fusion rules for $\V{+}$.
\bt{T3}{\rm(\cite{A})} Let $M$, $N$ and $L$ be irreducible
$\M{+}$-modules.

\noindent
{\rm(i)} If $M=\M{+}$, then $N_{{\M{+}}{ N}}^{L}=\delta_{N,L}$.

\noindent
{\rm(ii)} If $M=\M{-}$, then $N_{{\M{-}}{ N}}^{L}$ is $0$ or $1$, and
$N_{{\M{-}}{ N}}^{L}=1$ if and only if the pair $(N,L)$ is one of
 the following pairs:
\eqn
&&(\M{\pm},\M{\mp}),\ (\Mt{\pm},\Mt{\mp}),\\ 
&&(\Ml{\lm},\Ml{\mu})\hbox{ for $\lm,\mu\in\h-\{0\}$ such that 
$\B\lm,\lm\K=\B\mu,\mu\K$.}
\eeqn 

\noindent
{\rm(iii)} If $M=\Ml{\lm}$ for $\lm\in\h-\{0\}$, then $N_{{\Ml{\lm}}{ N}}^{L}$
is $0$ or $1$,  an
$N_{{\Ml{\lm}}{N}}^{L}=1$ if and only if the pair $(N,L)$ is one of
the following pairs:
\eqn
&&(\M{\pm},\Ml{\mu})\ (\Ml{\mu},\M{\pm})\hbox{ for $\mu\in\h-\{0\}$ 
such that $\B\lm,\lm\K=\B\mu,\mu\K$,}\\
&&(\Ml{\mu},\Ml{\nu})\hbox{ for $\mu,\nu\in\h-\{0\}$ such that 
$\B\nu,\nu\K=\B\lm\pm\mu,\lm\pm\mu\K$,}\\
&&(\Mt{\pm},\Mt{\pm}),\ (\Mt{\pm},\Mt{\mp}).
\eeqn

\noindent
{\rm(iv)} If $M=\Mt{+}$, then $N_{{\Mt{+}}{N}}^{L}$ is $0$ or $1$, and
$N_{{\Mt{+}}{ N}}^{L}=1$ if and only if the pair $(N,L)$ is one of the following pairs:
\eqn
&&(\M{\pm},\Mt{\pm}),\ (\Mt{\pm},\M{\pm}),\\ 
&&(\Ml{\lm},\Mt{\pm}),\ (\Mt{\pm},\Ml{\lm})\hbox{ for $\lm\in\h-\{0\}$.}
\eeqn

\noindent
{\rm(v)} If $M=\Mt{-}$, then $N_{{\Mt{-}} N}^{L}$ is $0$ or $1$, and
$N_{{\Mt{-}} N}^{L}=1$ if and only if the pair $(N,L)$ is one of the 
following pairs:
\eqn
&&(\M{\pm},\Mt{\mp}),\ (\Mt{\pm},\M{\mp}),\\ 
&&(\Ml{\lm},\Mt{\pm}),\ (\Mt{\pm},\Ml{\lm})\hbox{ for $\lm\in\h-\{0\}$.}
\eeqn
\et

Next we discuss the contragredient modules of irreducible $\V{+}$-modules. 
We shall prove the following proposition.
\bp{Con1}{\rm (i)} If $k$ is even, then all irreducible $\V{+}$-modules are
self-dual, that is , $W\cong W'$ as $\V{+}$-modules.

\noindent
{\rm (ii)} If $k$ is odd, then 
\eqn
(\Va{\pm})'\cong \Va{\mp}, (\V{T_{1},\pm})'\cong\V{T_{2},\pm},
(\V{T_{2},\pm})'\cong\V{T_{1},\pm}
\eeqn
 and others are self-dual.
\ep

To prove the proposition, we use Zhu's theory (see \cite{Z}). Let $V$ be a vertex 
operator algebra. The Zhu's algebra $A(V)$ associated with $V$ is a quotient space
of $V$ by the subspace $O(V)$ which is spanned by vectors of the form 
\eqn
a\circ b=\Res_{z}\frac{(1+z)^{\swt(a)}}{z^{2}} Y(a,z)b
\eeqn
for homogeneous element $a\in V$ and $b\in V$. The product of $A(V)$ is induced from the bilinear map $*:V\times V\to V$ which is defined by
\eqn
a*b=\Res_{z}\frac{(1+z)^{\swt(a)}}{z} Y(a,z)b
\eeqn
for homogeneous element $a\in V$ and $b\in V$. Let
$M$ be an irreducible $V$-module. Then there is a constant $h\in\C$ such that
$M$ has a direct sum decomposition $M=\ops_{n\in\N}M_{n}$, $M_{n}=\{\,v\in
M\,|\,L(0)v=(h+n)v\,\}$ for $n\in\N$. Then the action
$o(a)u=a_{\swt(a)-1}u$ for $a\in V$ and $u\in M$ induces an $A(V)$-module
structure on $M_{0}$ which is called the top level of $M$, and $M_{0}$ is
irreducible as $A(V)$-module. Furthermore if two irreducible $V$-modules $M$
and $N$ have top levels $M_{0}$ and
$N_{0}$ which are isomorphic to each other as $A(V)$-modules, then $M$ and
$N$ are isomorphic as $V$-module. 

Suppose that $k\neq1$. 
Then in \cite{DN2}, it is proved that the Zhu's algebra
$A(\V{+})$ is generated by three elements $[\w]$, $[J]$ and $[E]$, where
$[a]$ means the image $a+O(\V{+})$ in $ A(\V{+})$ of $a\in\V{+}$, and
$J=h(-1)^{4}\1-2h(-3)h(-1) \1+(3/2)h(-2)^{2}\1$ and
$E=e_{\al}+e_{-\al}$. Hence for an irreducible $\V{+}$-module $M$, to find the
irreducible module which is isomorphic to $M'$, it is enough to see the actions of
$[\w]$, $[J]$ and $[E]$ on the top level of $M'$. Since the top level of an irreducible
$\V{+}$-module is one dimensional, they act on the top level as scalar multiple. By
the construction of irreducible $\V{+}$-modules, we have the following table.

\begin{center}
\begin{tabular}{|c|c|c|c|c|c|}\hline
&$\V{+}$&$\V{-}$&$\Vl{\lm_{r}}$ $(1\leq r\leq k-1)$&$\Va{+}$&$\Va{-}$\\
\hline\hline 
$\w$&$0$&$1$&$r^{2}/4k$&$k/4$&$k/4$\\ 
$J$&$0$&$-6$&$(r^{2}/2k)^{2}-r^{2}/4k$&$k^{4}/4-k^{2}/4$&$k^{4}/4-k^{2}/4$\\

$E$&$0$&$0$&$0$&1&$-1$\\ \hline
\end{tabular}
\end{center}
\begin{center}
\begin{tabular}{|c|c|c|c|c|}\hline
&$\V{T_{1},+}$&$\V{T_{1},-}$&$\V{T_{2},+}$&$\V{T_{2},-}$\\
\hline\hline 
$\w$&$1/16$&$9/16$&$1/16$&$9/16$\\ 
$J$&$3/128$&$-45/128$&$3/128$&$-45/128$\\
$E$&$2^{-2k+1}$&$-2^{-2k+1}(4k-1)$&$-2^{-2k+1}$&$2^{-2k+1}(4k-1)$\\ \hline
\end{tabular}
\end{center}
\begin{center}
{{\bf Table 1.} Actions of  $\w$, $J$ and $E$ on the top level}
\end{center}

Now we prove Proposition \ref{Con1}.

{\it Proof of Proposition \ref{Con1}.} Firs we consider the case $k\neq1$. Let $W$
be an irreducible $\V{+}$-module. Set the top level $W_{0}=\C v$, and the top level
of the contragredient module $W'_{0}=\C v'$. By the definition of a contragredient
module, if $a\in\V{+}$ satisfies that $L(0)a=\wt(a)a$ and
$L(1)a=0$, we have 
$
\B o(a)v',v\K=(-1)^{\swt(a)}\B v',o(a)v\K,
$
and hence
\eqa
\B o(\w)v',v\K=\B v',o(\w)v\K, \B o(J)v',v\K=\B v',o(J)v\K,\B o(E)v',
v\K=(-1)^{k}\B v',o(E)v\K.\label{Act1}
\eeqa 
Therefore by Table $1$, we have Proposition \ref{Con1} for $k\neq1$. 

If $k=1$, then the dimension of the top level $(\V{-})_{0}$ is two
and others are one.  Hence we see that $(\V{-})'\cong\V{-}$ because the dimension
of $(\V{-})'_{0}$ is two.  Since for irreducible $\V{+}$-modules except $\V{-}$ 
Table $1$ is valid, we may apply same arguments of the case $k\neq1$ to such
irreducible modules. Therefore Table 1 and (\ref{Act1}) shows that
$$
(\V{+})'\cong\V{+},(\Va{\pm})'\cong \Va{\mp},(\V{T_{1},\pm})'\cong\V{T_{2},\pm},(\V{T_{2},\pm})'\cong\V{T_{1},\pm}.
$$
This proves Proposition \ref{Con1} for $k=1$.\qed 


\section{Fusion rules for $\V{+}$}\label{S4}
In Section \ref{S3.2}, we give 
irreducible decompositions of irreducible $\V{+}$-modules as $\M{+}$-modules,
and prove that every fusion rules for $\V{+}$ are zero or one with the help of fusion
rules for $\M{+}$. The main theorem is stated in Section \ref{S4.4}. The rest of
 sections is devoted to the proof of the theorem and it is divided into two
cases; one is the case that all modules are untwisted types (Section \ref{S4.1}) and
the other is the case that some irreducible module is twisted type (Section
\ref{S4.2}). 

\subsection{Irreducible decompositions of irreducible $\V{+}$-modules \\
as $\M{+}$-modules}\label{S3.2}

Since $\Vl{\lm}=\Ops _{m\in\Z}\Ml{\lm+m\al}$ for $\lm\in L^{\circ}$ and 
$\Ml{\mu}$ is irreducible for $\M{+}$ if $\mu\neq 0$, $\Vl{\lm_{r}}$ ($1\leq r\leq
k-1$) has an irreducible decompositions for $\M{+}$: 
\eqa \Vl{\lm_{r}}\cong\Ops _{m\in\Z}\Ml{\lm_{r}+m\al},\label{dec5}
\eeqa

For a nonzero $\lm\in\h$, we consider the subspace $(\M{+}\ots
(e_{\lm}\pm e_{-\lm}))\ops(\M{-}\ots (e_{\lm}\mp e_{-\lm}))$ on
$\Ml{\lm}\ops\Ml{-\lm}$. Since the action of $\M{+}$ of
$\Ml{\lm}\ops\Ml{-\lm}$ commutes the action of $\th$, the subspaces 
$(\M{+}\ots(e_{\lm}\pm e_{-\lm}))\ops(\M{-}\ots (e_{\lm}\mp
e_{-\lm}))$ are $\M{+}$-submodules. In fact we have the following
proposition.
\bl{P19} For a nonzero $\lm\in \h$, $\M{+}$-submodules 
$(\M{+}\ots (e_{\lm}\pm e_{-\lm}))\ops(\M{-}\ots (e_{\lm}\mp e_{-\lm}))$ of 
$\Ml{\lm}\ops\Ml{-\lm}$  are isomorphic to $\Ml{\lm}$. 
\el
\proof 
 Define a linear map $\phi_{\lm}$ by
\eqa 
\phi_{\lm}:(\M{+}\ots(e_{\lm}+e_{-\lm}))
\ops(\M{-}\ots (e_{\lm}-e_{-\lm}))&\to& \Ml{\lm}\label{ISO1}\\
u\ots(e_{\lm}+e_{-\lm})
+v\ots(e_{\lm}-e_{-\lm})&\mapsto& (u+v)\ots e_{\lm},\nn
\eeqa
for $u\in\M{+}$ and $v\in\M{-}$. Then the linear map
$\phi_{\lm}$ is an
 injective $\M{+}$-module homomorphism. Since $\Ml{\lm}$ is irreducible
for $\M{+}$, the homomorphism is in fact an isomorphism. Hence
$\M{+}\ots(e_{\lm}+e_{-\lm})\ops\M{-}\ots (e_{\lm}-e_{-\lm})$ is isomorphic to
$\Ml{\lm}$ as $\M{+}$-module.  We can also prove that
$\M{+}\ots(e_{\lm}-e_{-\lm})\ops\M{-}\ots (e_{\lm}+e_{-\lm})$ is isomorphic
to $\Ml{\lm}$ as $\M{+}$-module in the same way.\qed
\vs

We give irreducible decompositions of irreducible $\V{+}$-modules for $\M{+}$;
\bp{IRD1} Each irreducible $\V{+}$-modules decompose into direct sums of
irreducible $\M{+}$-modules as follows{\rm;}
\eqa  
\V{\pm}&\cong&\M{\pm}\ops\Ops_{m=1}^{\infty}\Ml{m\al},\label{dec1}\\
\Vl{\lm_{r}}&\cong&\Ops _{m\in\Z}\Ml{\lm_{r}+m\al}\hbox{ for $1\leq
r\leq k-1$},\label{dec2}\\ 
V_{\frac{\al}{2}+L}^{\pm}&\cong&
\Ops_{m=0}^{\infty}\Ml{{\frac{\al}{2}}+m\al},\label{dec3}\\
\V{T_{i},\pm}&\cong& \Mt{\pm}\hbox{  for $i=1,2$.}\label{dec4}
\eeqa
\ep
\proof Irreducible decompositions of $\Vl{\lm_{r}}$ ($1\leq r\leq k-1$) and
$\V{T_{i},\pm}$ ($i=1,2)$ have already given by (\ref{dec5}) and (\ref{ISO2})
respectively. We see that $\V{\pm}$ and $\Va{\pm}$ have direct sum
decompositions
\eqn  
\V{\pm}&=&\Ops_{m=0}^{\infty}((\M{+}\ots
(e_{m\al}\pm e_{-m\al}))\ops(\M{-}\ots (e_{m\al}\mp e_{-m\al}))),\nn\\
V_{\frac{\al}{2}+L}^{\pm}&=&\Ops_{m=0}^{\infty}((\M{+}\ots
(e_{{\frac{\al}{2}}+m\al}\pm e_{-{\frac{\al}{2}}-m\al}))\ops(\M{-}\ots
(e_{{\frac{\al}{2}}+m\al}\mp e_{-{\frac{\al}{2}}-m\al}))).\nn
\eeqn
Hence Lemma \ref{P19} shows that these direct sum decompositions give
irreducible decompositions of $\V{\pm}$ and $\Va{\pm}$.\qed\vs
 
By Proposition \ref{IRD1}, one see that for any irreducible $\V{+}$-module $W$, the
multiplicity of an irreducible $\M{+}$-module in $W$ is at most one.

Using these irreducible decompositions (\ref{dec1})-(\ref{dec4}),
Theorem \ref{T3} and Corollary \ref{C8}, we can show that all fusion rules for
$\V{+}$ are at most one:
\bp{P49}Let $W^{1},W^{2}$ and $W^{3}$ be irreducible $\V{+}$-module.

\noindent
{\rm(1)} The fusion rule $N_{W^{1}W^{2}}^{W^{3}}$ is zero or one.

\noindent
{\rm(2)} If all $W^{i}$ $(i=1,2,3)$ are twisted type modules, then the fusion rule
$N_{W^{1}W^{2}}^{W^{3}}$ is zero.

\noindent
{\rm(3)} If one of $W^{i}$ $(i=1,2,3)$ is twisted type module and others are of
untwisted types, then the fusion rule $N_{W^{1}W^{2}}^{W^{3}}$ is zero.
\ep
\proof 
Suppose that $W^{1}$ and $W^{2}$ have irreducible $\M{+}$-submodules $M$
 and $N$ respectively and that $W^{3}$ has an irreducible decomposition
$W^{3}=\ops_{i} M^{i}$ as $\M{+}$-module. By (\ref{Ineq2}), we have an
inequality
\eqa
\dim I_{\V{+}}\fusion{W^{1}}{W^{2}}{W^{3}}\leq \sum_{i}\dim
I_{\M{+}}\fusion{M}{N}{M^{i}}.\label{Ineq1}
\eeqa

If $W^{1}, W^{2}$ and $W^{3}$ are of twisted type or if $W^{1}$ is of twisted type
and  $W^{2}$ and $W^{3}$ are of untwisted type, then by Theorem \ref{T3} (iv),
(v) and (\ref{dec1})-(\ref{dec4}), we see that the fusion rule for
$\M{+}$ of type $\fusion{M}{N}{M^{i}}$ is zero for any $i$. Hence (\ref{Ineq1})
implies that the fusion rule $N_{W^{1}W^{2}}^{W^{3}}$ is zero. Since the
contragredient module of an (un)twisted type module is of (un)twisted type, (2)
and (3) follows from Proposition \ref{P9}.

By (2),(3) and Proposition \ref{P9}, to show (1), it suffices to prove that if
$W^{1}$ is untwisted type module and both $W^{2}$ and $W^{3}$ are of twisted
types or of untwisted types, then the fusion rule $N_{W^{1}W^{2}}^{W^{3}}$ is zero
or one. 

If $W^{1}$ is untwisted type module and $W^{2}$ and $W^{3}$ are of twisted types,
then by Theorem \ref{T3} (i)-(iii) and irreducible decompositions
(\ref{dec1})-(\ref{dec4}), we see that the fusion rule for $\M{+}$ of type
$\fusion{M}{W^{2}}{W^{3}}$ is zero or one for any irreducible
$\M{+}$-submodules $M$ of $W^{1}$. Hence Corollary \ref{C8} shows that the
fusion rule
$N_{W^{1}W^{2}}^{W^{3}}$ is zero or one.

Now we turn to the case all $W^{i}$ $(i=1,2,3)$ are of untwisted types. We consider
the following three cases separately; (i)
$W^{1}=\V{\pm}$, (ii) $W^{1}=\Va{\pm}$ and (iii) $W^{1}=\Vl{\lm_{r}}$ for $1\leq
r\leq k-1$. Let $W^{3}=\ops_{i}M^{i}$ be the irreducible decomposition of $W^{3}$
for $\M{+}$. Then it suffices to prove that the right-hand side of (\ref{Ineq1}) is at
most one for some $\M{+}$-submodules
$M$ of $W^{1}$ and $N$ of $W^{2}$. 

(i) $W^{1}=\V{\pm}$ cases: Take $M=\M{\pm}$. By
(\ref{dec1})-(\ref{dec3}), we can take $N$ to be isomorphic to $\Ml{\lm}$ for
some
$\lm\in L^{\circ}$.  Then by Theorem \ref{T3} (i) and (ii), the fusion rule for
$\M{+}$ of type $\fusion{M}{N}{M^{i}}$ is one if and only if $M^{i}$ is
isomorphic to $\Ml{\lm}$.  Since the multiplicity of $\Ml{\lm}$ in the
irreducible decomposition of $W^{3}$ is at most one, we see that the
right-hand side of (\ref{Ineq1}) is zero or one.  

(ii) $W^{1}=\Va{\pm}$ case: Take $M\cong\Ml{\al/2}$. If $N$ is isomorphic to
$\Ml{\lm}$ for some $\lm\in L^{\circ}$, then by Theorem
\ref{T3} (iii), we see that the fusion rule for
$\M{+}$ of type $\fusion{M}{N}{M^{i}}$ is one if and only if $M^{i}$ is
isomorphic to $\Ml{\lm+\al/2}$ or $\Ml{\lm-\al/2}$. If
$W^{2}=\Va{\pm}$, then by taking $\lm=\al/2$, we see that the right-hand side
of (\ref{Ineq1}) is zero unless $W^{3}$ is $\V{+}$ or $\V{-}$. So these cases  and
the cases $W^{2}=\V{\pm}$ reduce to the case (i) by means of
Proposition
\ref{P9}. Therefore to prove (1) in the case $W^{1}=\Va{\pm}$,  it is enough to
consider the case $W^{2}=\Vl{\lm_{r}}$ for some $1\leq r\leq k-1$. Then by taking
$\lm=\lm_{r}$, we see that the right-hand side of (\ref{Ineq1}) is zero unless
$W^{3}$ is
$\Vl{\lm_{k-r}}$. By Corollary \ref{C8} and Proposition \ref{Con1}, the fusion rules
of types
$\fusion{\Va{\pm}}{\Vl{\lm_{r}}}{\Vl{\lm_{k-r}}}$ is equal those of types
$\fusion{\Vl{\lm_{r}}}{\Vl{\lm_{k-r}}}{(\Va{\pm})'}$ respectively. Hence we have
to show that the right-hand side of (\ref{Ineq1}) is at most one when
$W^{1}=\Vl{\lm_{r}}, W^{2}=\Vl{\lm_{k-r}}$ and $W^{3}=(\Va{\pm})'$. We take
$M=\Ml{\lm_{r}}$ and $N=\Ml{\lm_{k-r}}$. Since by Theorem \ref{T3} the fusion
rule for $\M{+}$ of type  $\fusion{M}{N}{\Ml{\al/2+m\al}}$ is $\del_{m,0}$ for
$m\in\N$, (\ref{dec3}) shows that the the right-hand side of (\ref{Ineq1}) is
at most one.  

(iii) $W^{1}=\Vl{\lm_{r}}$ case for $1\leq r\leq k-1$: 
By Proposition \ref{P9} and results of (i) and (ii),  to prove (1) in the case, it is
sufficient to consider the case $W^{2}=\Vl{\lm_{s}}$ for $1\leq s\leq k-1$. Then we
can take
$M=\Ml{\lm_{r}}$ and $N=\Ml{\lm_{s}}$. Hence by
Theorem \ref{T3} (iii), we see that the fusion rule for $\M{+}$ of type
$\fusion{M}{N}{M^{i}}$ is one if and only if $M^{i}$ is isomorphic to
$\Ml{\lm_{r}+\lm_{s}}$ or $\Ml{\lm_{r}-\lm_{s}}$. By
(\ref{dec1})-(\ref{dec3}), one see that for $\mu,\nu\in L^{\circ}$
$\Ml{\mu}$ and $\Ml{\nu}$ have multiplicity one in $W^{3}$, then
$\mu+\nu\in L$ or $\mu-\nu\in L$. But
$(\lm_{r}+\lm_{s})+(\lm_{r}-\lm_{s})$ and
$(\lm_{r}+\lm_{s})-(\lm_{r}-\lm_{s})$ are not in $L$. Hence
by (\ref{dec1})-(\ref{dec3}), we see that the right-hand side  of (\ref{Ineq1}) is
zero or one.\qed
\vs


\subsection{Main Theorem}\label{S4.4}
Here we give the main theorem. The proof is given in Section \ref{S4.1} and
\ref{S4.2}:
\bt{T4} Let $W^{1},W^{2}$ and $W^{3}$ be irreducible $\V{+}$-modules. Then
{\rm(1)} the fusion rule
$N_{W^{1}W^{2}}^{W^{3}}$ is zero or one and {\rm(2)} the fusion rule
$N_{W^{1}W^{2}}^{W^{3}}$ is one if and only if $W^{i}$ $(i=1,2,3)$ satisfy following
cases{\rm :} 
 
\noindent
{\rm(i)} $W^{1}=\V{+}$ and $W^{2}\cong W^{3}$.

\noindent
{\rm(ii)} $W^{1}=\V{-}$ and the pair $(W^{2},W^{3})$ is one of the pairs
\eqn&&(\V{\pm},\V{\mp}),\ (V_{\al/2+L}^{\pm},V_{\al/2+L}^{\mp}),\
(\V{T_{1},\pm},\V{T_{1},\mp}),\ (\V{T_{2},\pm},\V{T_{2},\mp}),\\
&&(\Vl{\lm_{r}},\Vl{\lm_{r}})\ for\ 1\leq
r\leq k-1.
\eeqn

\noindent
{\rm(iii)} $W^{1}=\Va{+}$ and the pair $(W^{2},W^{3})$ is one of the pairs
\eqn
&&(\V{\pm},\Va{\pm}),\  ((\Va{\pm})',\V{\pm}),\ 
((\V{T_{1},\pm})',\V{T_{1},\pm}),\ ((\V{T_{2},\pm})',\V{T_{2},\mp}),\\ 
&&(\Vl{\lm_{r}},\Vl{\al/2-\lm_{r}})\ for\
1\leq r\leq k-1.
\eeqn

\noindent
{\rm(iv)} $W^{1}=\Va{-}$ and the pair $(W^{2},W^{3})$ is one of the pairs
\eqn
&&(\V{\pm},\Va{\mp}),\  ((\Va{\pm})',\V{\mp}),\ 
((\V{T_{1},\pm})',\V{T_{1},\mp}),\ ((\V{T_{2},\pm})',\V{T_{2},\pm}),\\ 
&&(\Vl{\lm_{r}},\Vl{\al/2-\lm_{r}})\ for\
1\leq r\leq k-1.
\eeqn

\noindent
{\rm(v)} $W^{1}=\Vl{\lm_{r}}$ for $1\leq r\leq k-1$
 and the pair $(W^{2},W^{3})$ is one of the pairs
\eqn
&&(\V{\pm},\Vl{\lm_{r}}),\ (\V{\pm},\Vl{\lm_{r}}),\ 
(\Va{\pm},\Vl{\al/2-\lm_{r}}),\ (\Vl{\al/2-\lm_{r}},\Va{\pm}),\\ 
&&(\Vl{\lm_{s}},\Vl{\lm_{r}\pm\lm_{s}})\hbox{ for $1\leq s\leq k-1$,}\\
&&(\V{T_{1},\pm},\V{T_{1},\pm}),\ (\V{T_{1},\pm},\V{T_{1},\mp}),\ 
(\V{T_{2},\pm},\V{T_{2},\pm}),\ (\V{T_{2},\pm},\V{T_{2},\mp})\hbox{ if
$r$ is even},\\  
&&(\V{T_{1},\pm},\V{T_{2},\pm}),\
(\V{T_{1},\pm},\V{T_{2},\mp}),\ (\V{T_{2},\pm},\V{T_{1},\pm}),\
(\V{T_{2},\pm},\V{T_{1},\mp})\hbox{ if $r$ is odd}.
\eeqn

\noindent
{\rm(vi)} $W^{1}=(\V{T_{1},+})'$ and the pair $(W^{2},W^{3})$ is one of the pairs
\eqn &&(\V{\pm},(\V{T_{1},\pm})'),\ (\V{T_{1},\pm},\V{\pm}),\ 
(\Va{\pm},\V{T_{1},\pm}),\ ((\V{T_{1},\pm})',(\Va{\pm})'),\\ 
&&(\Vl{\lm_{r}}, (\V{T_{1},\pm})')\hbox{ and }(\V{T_{1},\pm},\Vl{\lm_{r}})\hbox{
for $1\leq r\leq k-1$ and $r$ is even},\\ 
&&(\Vl{\lm_{r}},(\V{T_{2},\pm})')\hbox{ and
}(\V{T_{2},\pm},\Vl{\lm_{r}})\hbox{ for $1\leq
r\leq k-1$ and $r$ is odd}.
\eeqn

\noindent
{\rm(vii)} $W^{1}=(\V{T_{1},-})'$ and the pair $(W^{2},W^{3})$ is one of the pairs
\eqn &&(\V{\pm},(\V{T_{1},\mp})'),\ (\V{T_{1},\pm},\V{\mp}),\ 
(\Va{\pm},\V{T_{1},\mp}),\ ((\V{T_{1},\pm})',(\Va{\mp})'),\\ 
&&(\Vl{\lm_{r}}, (\V{T_{1},\pm})')\hbox{ and }(\V{T_{1},\pm},\Vl{\lm_{r}})\hbox{
for $1\leq r\leq k-1$ and $r$ is even},\\ 
&&(\Vl{\lm_{r}},(\V{T_{2},\pm})')\hbox{ and
}(\V{T_{2},\pm},\Vl{\lm_{r}})\hbox{ for $1\leq
r\leq k-1$ and $r$ is odd}.
\eeqn

\noindent
{\rm(viii)} $W^{1}=(\V{T_{2},+})'$ and the pair $(W^{2},W^{3})$ is one of the pairs
\eqn &&(\V{\pm},(\V{T_{2},\pm})'),\ (\V{T_{2},\pm},\V{\pm}),\ 
(\Va{\pm},\V{T_{2},\mp}),\ ((\V{T_{2},\pm})',(\Va{\mp})'),\\ 
&&(\Vl{\lm_{r}}, (\V{T_{2},\pm})')\hbox{ and }(\V{T_{2},\pm},\Vl{\lm_{r}})\hbox{
for $1\leq r\leq k-1$ and $r$ is even},\\ 
&&(\Vl{\lm_{r}},(\V{T_{1},\pm})')\hbox{ and
}(\V{T_{1},\pm},\Vl{\lm_{r}})\hbox{ for $1\leq
r\leq k-1$ and $r$ is odd}.
\eeqn

\noindent
{\rm(ix)} $W^{1}=(\V{T_{2},-})'$ and the pair $(W^{2},W^{3})$ is one of the pairs
\eqn &&(\V{\pm},(\V{T_{2},\mp})'),\ (\V{T_{2},\pm},\V{\mp}),\ 
(\Va{\pm},\V{T_{2},\pm}),\ ((\V{T_{2},\pm})',(\Va{\pm})'),\\ 
&&(\Vl{\lm_{r}}, (\V{T_{2},\pm})')\hbox{ and }(\V{T_{2},\pm},\Vl{\lm_{r}})\hbox{
for $1\leq r\leq k-1$ and $r$ is even},\\ 
&&(\Vl{\lm_{r}},(\V{T_{1},\pm})')\hbox{ and
}(\V{T_{1},\pm},\Vl{\lm_{r}})\hbox{ for $1\leq
r\leq k-1$ and $r$ is odd}.
\eeqn
\et

Since (1) of the main theorem has already proved in Proposition \ref{P49}, to
prove the theorem, it is enough to show that for irreducible $\V{+}$-modules
$W^{1},W^{2}$ and
$W^{3}$,  the fusion rule $N_{W^{1}W^{2}}^{W^{3}}$ is nonzero if and only if the
triple $(W^{1},W^{2},W^{3})$ satisfy indicated cases in the theorem.  In Section
\ref{S4.1}, we prove this in the case all $W^{i}$ $(i=1,2,3)$ are untwisted type
modules, and in Section \ref{S4.2} we do in the case one of $W^{i}$ $(i=1,2,3)$ is
twisted type modules.  To show "if" part, we shall construct nonzero intertwining
operators explicitly.
\subsection{Fusion rules for untwisted type modules}\label{S4.1}
We construct nonzero intertwining operators for untwisted type modules.
For this purpose, we review intertwining operators for $V_{L}$ following \cite{DL}.

Let $\lm,\mu\in L^{\circ}$. An intertwining
operator of type $\fusion{\Vl{\lm}}{\Vl{\mu}}{\Vl{\lm+\mu}}$ is
constructed as follows:  As shown in Chapter 8 of \cite{FLM},  the
operator $\mathcal{Y}^{\circ}$ satisfies Jacobi identity and $L(-1)$-derivative
property on
$V_{L^{\circ}}$ for
$\beta\in L$, $\lm\in L^{\circ},a\in\Ml{\beta}$ and $u\in\Ml{\lm}$:
\eqn
&&z_{0}^{-1}\delta\left(
{\frac{z_{1}-z_{2}}{z_{0}}}\right)Y(a,z_{1})\mathcal{Y}^{\circ}(u,z_{2})
-(-1)^{\B\beta,\lm\K}z_{0}^{-1}\delta\left({\frac{z_{2}-z_{1}}{
-z_{0}}}\right)\mathcal{Y}^{\circ}(u,z_{2})Y(a,z_{1})\\
&&{ }=z_{2}^{-1}\delta\left( {\frac{z_{1}-z_{0}}{z_{2}}}\right)
\mathcal{Y}^{\circ}(Y(a,z_{0})u,z_{2}),
\eeqn
\eqn
\frac{d}{dz}\mathcal{Y}^{\circ}(u,z)=\mathcal{Y}^{\circ}(L(-1)u,z).
\eeqn
Let $\pi_{\lm}$ $(\lm\in L^{\circ})$ be the linear endomorphism of
$V_{L^{\circ}}$ defined by $\pi_{\lm}(v)=e^{\B \lm,\mu\K \pi i}\,v$ for
$\mu\in L^{\circ}$ and $v\in\Ml{\mu}$.
Set $\mathcal{Y}_{r,s}(u,z)=\mathcal{Y}^{\circ}(u,z)
\pi_{\lm_{r}}|_{\Vl{\lm_{s}}}$ for $r,s\in\Z$ and $u\in\Vl{\lm_{r}}$.
Then the operator $\mathcal{Y}_{r,s}$ gives a nonzero intertwining operator for
$\V{}$ of type $\fusion{\Vl{\lm_{r}}}{\Vl{\lm_{s}}}{\Vl{\lm_{r}+\lm_{s}}}$ (see
\cite{DL}). 

\bp{P5} Fusion rules of the following types are nonzero{\rm ;}
\begin{enumerate}
\renewcommand{\labelenumi}{\rm (\roman{enumi})\ \ }
\item
$\fusion{\Vl{\lm_{r}}}{\Vl{\lm_{s}}}{\Vl{(\lm_{r}\pm\lm_{s})}}\hbox{ for
}1\leq r,s\leq k-1,$ 
\item$\fusion{\V{+}}{\V{\pm}}{\V{\pm}},
\fusion{\V{-}}{\V{\pm}}{\V{\mp}}$
and 
$\fusion{\V{\pm}}{\Vl{\lm_{r}}}{\Vl{\lm_{r}}}\hbox{ for }0\leq r\leq k-1,$
\item
$\fusion{\V{+}}{\Va{\pm}}{\Va{\pm}},
\fusion{\V{-}}{\Va{\pm}}{\Va{\mp}}$
and $\fusion{\Va{\pm}}{\Vl{\lm_{r}}}{\Vl{(\al/2-
\lm_{r})}}\hbox{ for }0\leq r\leq k-1.$ 
\end{enumerate}
\ep
\proof Since $(\V{},Y),(\Va{},Y)$ and $(\Vl{\lm_{r}},Y)$ $(1\leq r\leq k-1$) are
irreducible $\V{}$-modules, the vertex operator $Y$ gives nonzero intertwining
operators for
$\V{}$ of types $$\fusion{\V{}}{\V{}}{\V{}},
\fusion{\V{}}{\Va{}}{\Va{}}\hbox{ and }\fusion{\V{}}{\Vl{\lm_{r}}}{\Vl{\lm_{r}}}$$ 
for $1\leq r\leq k-1$.
Hence $Y(a,z)u$ is nonzero for any nonzero $a\in\V{}$ and nonzero
$u\in\Vl{\lm_{s}}$ ($s\in\Z$) by Corollary \ref{C8}. Therefore since $\th
Y(a,z)\th=Y(\th(a),z)$ for $a\in\V{}$, 
$Y$ gives nonzero intertwining operators for $\V{+}$ of types
$$\fusion{\V{+}}{\V{\pm}}{\V{\pm}},
\fusion{\V{-}}{\V{\pm}}{\V{\mp}},\fusion{\V{+}}{\Va{\pm}}{\Va{\pm}},\fusion{\V{-}}{\Va{\pm}}{\Va{\mp}}\hbox{ and }\fusion{\V{\pm}}{\Vl{\lm_{r}}}{\Vl{\lm_{r}}}$$
for $1\leq r\leq k-1.$

 Next we show that fusion rules of types
$\fusion{\Vl{\lm_{r}}}{\Vl{\lm_{s}}}{\Vl{(\lm_{r}\pm\lm_{s})}}$
for $r,s\in\Z$ are nonzero.
Define $\mathcal{Y}_{r,-s}\circ\th$ by
$(\mathcal{Y}_{r,-s}\circ\th)(u,z)v=\mathcal{Y}_{r,-s}(u,z)\th(v)$ for
$u\in\Vl{\lm_{r}}$ and $v\in\Vl{\lm_{s}}$. 
Then $\mathcal{Y}_{r,s}$ is a nonzero intertwining operator for $\V{+}$ of type
$\fusion{\Vl{\lm_{r}}}{\Vl{\lm_{s}}}{\Vl{(\lm_{r}-\lm_{s})}}$ since $\th$ commutes the action of $\V{+}$. 
This proves that fusion rules of types
$\fusion{\Vl{\lm_{r}}}{\Vl{\lm_{s}}}{\Vl{(\lm_{r}\pm\lm_{s})}}$ are nonzero
for any $r,s\in\Z$. 

Finally we show that fusion rule of type
$\fusion{\Va{\pm}}{\Vl{\lm_{r}}}{\Vl{(\al/2-\lm_{r})}}$
for $r\in\Z$ are nonzero.
Since $\Va{}$ is an irreducible $\V{}$-module for any
$r\in\Z$, Corollary \ref{C8} shows that $\mathcal{Y}_{k,r}(u,z)v$ is nonzero for
any nonzero $u\in\Va{}$ and nonzero $v\in\Vl{\lm_{r}}$. 
Hence $(\mathcal{Y}_{k,-r}\circ\th)(u,z)v$ is also nonzero for any nonzero
$u\in\Va{}$ and nonzero $v\in\Vl{\lm_{r}}$. Therefore $\mathcal{Y}_{k,-r}\circ\th$ gives nonzero
intertwining operators of types $\fusion{\Va{\pm}}{\Vl{\lm_{r}}}{\Vl{(\al/2-\lm_{r})}}$.\qed
\vs

Next we show that for untwisted type modules $W^{i}$ ($i=1,2,3$), if the fusion
rule $N_{W^{1}W^{2}}^{W^{3}}$ is nonzero, then the type
$\fusion{W^{1}}{W^{2}}{W^{3}}$ is given from types in Proposition \ref{P5} by using
Proposition \ref{P1}. For this it suffices to prove the following proposition.
\bp{P31} Let $W^{1},W^{2}$ and $W^{3}$ be untwisted type modules. 
Then the fusion rule $N_{W^{1}W^{2}}^{W^{3}}$ is zero if $W^{i}$ $(i=1,2,3)$ satisfy
following cases{\rm:}
\begin{enumerate}
\renewcommand{\labelenumi}{\rm (\roman{enumi})\ \ }
\item $W^{1}=\V{+}$, and $W^{2},W^{3}$ are inequivalent.

\item $W^{1}=\V{-}$ and the pair $(W^{1},W^{2})$ is one of following pairs
\eqn
(W^{2},W^{3})&=&(\V{-},\V{-}),\
(\Va{\pm},\Va{\pm}),\\
&&(\Vl{\lm_{r}},\Vl{\lm_{s}})\hbox{ for $1\leq r,s\leq
k-1$ and 
$\lm_{r}\neq\lm_{s}$},\\ 
&&(\V{-},\Vl{\lm_{r}}),\ (\Va{\pm},\Vl{\lm_{r}})\hbox{
for $1\leq r\leq k-1$}.
\eeqn

\item$W^{3}=\Va{\pm}$ and the pair $(W^{1},W^{2})$ is one of following pairs
\eqn
(W^{1},W^{2})&=&(\Vl{\lm_{r}},\Vl{\lm_{s}})\hbox{ for $1\leq r,s\leq k-1$
and $\lm_{r}+\lm_{s}\neq\al/2$},\\  
&&(\Va{\pm},\Vl{\lm_{r}})\hbox{ and $1\leq r\leq k-1$}.
\eeqn

\item $W^{1}=\Vl{\lm_{r}}$ for $1\leq r\leq k-1$ and the pair $(W^{1},W^{2})$ is
one of following pairs
\eqn
&&(W^{1},W^{2})=(\Vl{\lm_{s}},\Vl{\lm_{t}})\hbox{ for $1\leq s,t\leq
k-1$ and $\lm_{t}\neq\lm_{r}\pm\lm_{s}$ and $\lm_{s}-\lm_{r}$}. 
\eeqn
\end{enumerate}
\ep 
\proof Lemma \ref{Le1} proves the proposition in the case (i). 

Next we consider the cases (ii) except the pair $(\Va{\pm},\Va{\pm})$, (iii) and (iv). 
Let $W^{3}=\ops_{i}M^{i}$ be the irreducible decomposition of $W^{3}$ for $\M{+}$. 
Then we can find irreducible $\M{+}$-submodules $M$ of $W^{1}$ and $N$ of $W^{3}$ such that the fusion rule for $\M{+}$ of type $\fusion{M}{N}{M^{i}}$ is zero. Hence (\ref{Ineq1}) implies that
the fusion rule $N_{W^{1}W^{2}}^{W^{3}}$ is zero; for example, in the case $W^{1}=W^{2}=W^{3}=\V{-}$, we take $M=N=\M{-}$, etc..

It remains to prove the proposition in the cases $(W^{1},W^{2})=(\Va{\pm},\Va{\pm})$ of (ii). We prove that the
fusion rule of type $\fusion{\V{-}}{\Va{+}}{\Va{+}}$ is zero.
The case of type $\fusion{\V{-}}{\Va{-}}{\Va{-}}$ can be also proved in the similar way. 

Set
\eqa
\Va{+}[m]&=&\M{+}\ots(e_{\frac{\al}{2}+m\al}+e_{-(\frac{\al}{2}+m\al)})
\ops\M{-}\ots(e_{\frac{\al}{2}+m\al}-e_{-(\frac{\al}{2}+m\al)}),\label{VA1}
\eeqa 
for $m\in\N$. Note that $\Va{+}[m]$ is isomorphic to $\Ml{\al/2+m\al}$ as
$\M{+}$-module by Proposition \ref{P19}. Let
$\mathcal{Y}$ be an intertwining operator of type
$\fusion{\V{-}}{\Va{+}}{\Va{+}}$. By Theorem \ref{T3} (ii), we have
$\mathcal{Y}(u,z)v\in \Va{+}[0]((z))$ for $u\in\M{-}$ and $v\in\Va{+}[0]$. 
Let $\phi_{{\al}/{2}}:\Va{+}[0]\to\Ml{{\al}/{2}}$ be the $\M{+}$-module isomorphism defined in (\ref{ISO1}). For simplicity, we denote
$\phi=\phi_{\al/2}$. Then the operator 
$\phi\circ\mathcal{Y}\circ \phi^{-1}$ defined by
$(\phi\circ\mathcal{Y}\circ\phi^{-1})(u,z)v=\phi\mathcal{Y}(u,z)\phi^{-1}(v)$ for
$u\in\M{-}$ and $v\in\Ml{\al/2}$ gives an intertwining operator of type
$\fusion{\M{-}}{\Ml{\al/2}}{\Ml{\al/2}}$. Since the
dimension of 
$I_{{\M{+}}}\fusion{\M{-}}{\Ml{\al/2}}{\Ml{\al/2}}$
is one and the corresponding intertwining operator is given by a scalar multiple of
the vertex operator $Y$ of the $\M{}$-module $(\Ml{\al/2},Y)$, there exists a 
constant $d\in\C$ such that 
$$\mathcal{Y}(u,z)v=d\,\phi^{-1}Y(u,z)\phi(v)$$
for every $u\in\M{-}$ and $v\in\Va{+}[0]$. 
We write $\mathcal{Y}(u,z)=\sum_{n\in\Z}\Tilde{u}(n)z^{-n-1}$ $\Tilde{u}\in\End{\Va{+}}$ for $u\in\V{-}$.
Take $u=h(-1)\1$ and $v=e_{\al/2}+e_{-\al/2}$, then we have 
\eqa \tilde{h}(0)(e_{\frac{\al}{2}}+e_{-\frac{\al}{2}})&=&d \B h,\frac{\al}{2}\K
(e_{\frac{\al}{2}}+e_{-\frac{\al}{2}}),\label{E6}\\
\tilde{h}(-1)(e_{\frac{\al}{2}}+e_{-\frac{\al}{2}})&=&d
(h(-1)e_{\frac{\al}{2}}-h(-1)e_{-\frac{\al}{2}}),\label{E7}
\eeqa
where we denote $(\widetilde{h(-1)\1})(n)$ by $\tilde{h}(n)$ for $n\in\Z$. By
direct culculations, we see that 
\eqa
&&E_{k-1}(e_{\frac{\al}{2}}+e_{-\frac{\al}{2}})
=(e_{\frac{\al}{2}}+e_{-\frac{\al}{2}}),\label{E8}\\
&&E_{k}(h(-1)e_{\frac{\al}{2}}-h(-1)e_{-\frac{\al}{2}})
=\B h,\al\K(e_{\frac{\al}{2}}+e_{-\frac{\al}{2}}),\label{E9}
\eeqa
where $E=e_{\al}+e_{-\al}\in\V{+}$. Let $F=e_{\al}-e_{-\al}\in\V{-}$. Then by
Jacobi identity, we have a commutation relation
\eqa [E_{m},\tilde{h}(n)]=-\B h,\al\K \tilde{F}(m+n)\label{E10}
\eeqa
for $m,n\in\Z$. 
Hence (\ref{E6}) and (\ref{E8}) imply that $\tilde{F}(k-1)
(e_{{\al}/{2}}+e_{-{\al}/{2}})=0$ (take $m=k-1,n=0$ in (\ref{E10})). On the other hand, by
(\ref{E7}) and (\ref{E9}) we have 
\eqn
-\B h,\al\K \tilde{F}(k-1)(e_{\frac{\al}{2}}+e_{-\frac{\al}{2}})&=&
[E_{k},\tilde{h}(-1)](e_{\frac{\al}{2}}+e_{-\frac{\al}{2}})\\
&=&d \B h,\al\K(e_{\frac{\al}{2}}+e_{-\frac{\al}{2}}),
\eeqn
(take $m=k,n=-1$ in
(\ref{E10})). Therefore $d=0$. This implies that
$\mathcal{Y}(h(-1)\1,z)(e_{{\al}/{2}}+e_{-{\al}/{2}})=0$, and then Lemma
\ref{P1} shows $\mathcal{Y}=0$. Thus the fusion rule of type
$\fusion{\V{-}}{\Va{+}}{\Va{+}}$ is zero.\qed
\vs

Consequently, by Proposition \ref{P9}, Proposition \ref{P31}, Proposition
\ref{Con1}, Proposition \ref{P5} and Proposition \ref{P49}, we can determine
fusion rules for untwisted type modules.

\bp{MP1} Let $W^{1},W^{2}$ and $W^{3}$ be untwisted type $\V{+}$-modules. Then the
fusion rule
$N_{W^{1}W^{2}}^{W^{3}}$ is zero or one. The fusion rule $N_{W^{1}W^{2}}^{W^{3}}$ is
one if and only if $W^{i}$ $(i=1,2,3)$ satisfy the following cases{\rm :} 
 
\noindent
{\rm(i)} $W^{1}=\V{+}$ and $W^{2}\cong W^{3}$.

\noindent
{\rm(ii)} $W^{1}=\V{-}$ and the pair $(W^{2},W^{3})$ is one of pairs
\eqn
(\V{\pm},\V{\mp}),\ (V_{\al/2+L}^{\pm},V_{\al/2+L}^{\mp}),\ (\Vl{\lm_{r}},\Vl{\lm_{r}})\ for\ 1\leq
r\leq k-1.
\eeqn

\noindent
{\rm(iii)} $W^{1}=\Va{+}$ and the pair $(W^{2},W^{3})$ is one of pairs
\eqn
(\V{\pm},\Va{\pm}),\  ((\Va{\pm})',\V{\pm}),\ (\Vl{\lm_{r}},\Vl{\al/2-\lm_{r}})\ for\
1\leq r\leq k-1.
\eeqn

\noindent
{\rm(iv)} $W^{1}=\Va{-}$ and the pair $(W^{2},W^{3})$ is one of pairs
\eqn
(\V{\pm},\Va{\mp}),\  ((\Va{\pm})',\V{\mp}),\ (\Vl{\lm_{r}},\Vl{\al/2-\lm_{r}})\ for\
1\leq r\leq k-1.
\eeqn

\noindent
{\rm(v)} $W^{1}=\Vl{\lm_{r}}$ for $1\leq r\leq k-1$
 and the pair $(W^{2},W^{3})$ is one of pairs
\eqn
&&(\V{\pm},\V{\pm}),\ 
(\Va{\pm},\Vl{\al/2-\lm_{r}}),\ (\Vl{\al/2-\lm_{r}},\Va{\pm}),\\ 
&&(\Vl{\lm_{s}},\Vl{\lm_{r}\pm\lm_{s}})\hbox{ for $1\leq s\leq k-1$.}
\eeqn
\ep

\subsection{Fusion rules involving twisted type modules}\label{S4.2}
Set $\Pa=L^{\circ}\times \{1,2\}\times\{1,2\}$. 
We call $(\lm,i,j)\in \Pa$ a {\it quasi-admissible triple} if $\lm,i$ and $j$ satisfies
$$(-1)^{\B\lm,\al\K+\del_{i,j}+1}=1.$$
We denote the set of all quasi-admissible triples by $\Qa$. For a quasi-admissible
triple $(\lm,i,j)\in \Qa$,  we first construct an intertwining operator for
$\V{+}$ of type $\fusion{\Vl{\lm}}{\V{T_{i}}}{\V{T_{j}}}.$  

As shown in Chapter 9 of \cite{FLM} the operator $\mathcal{Y}^{\th}$ satisfies
twisted Jacobi identity and $L(-1)$-derivative property 
\eqa
&& z_{0}^{-1}\delta\left(
{\frac{z_{1}-z_{2}}{z_{0}}}\right)\mathcal{Y}^{\th}(a,z_{1})
\mathcal{Y}^{\th}(u,z_{2})
-(-1)^{\B\beta,\lm\K}z_{0}^{-1}\delta\left(
{\frac{z_{2}-z_{1}}{-z_{0}}}\right)\mathcal{Y}^{\th}(u,z_{2})
\mathcal{Y}^{\th}(a,z_{1})\nn\\
&&{ }=\frac{1}{2}\sum_{p=0,1}z_{2}^{-1}\delta\left((-1)^{p}
{\frac{(z_{1}-z_{0})^{1/2}}{z_{2}^{1/2}}}
\right)\mathcal{Y}^{\th}(Y(\th^{p}(a),z_{0})u,z_{2}),\label{Ja1}
\eeqa
and
\eqa
\frac{d}{dz}\mathcal{Y}^{\th}(u,z)=\mathcal{Y}^{\th}(L(-1)u,z)\label{Ld1}
\eeqa
for $\beta\in L,\lm\in L^{\circ}, a\in\Ml{\beta}$ and $u\in\Ml{\lm}$. 
Then we have following lemma.
\bl{L3.4}{\rm(1)} The intertwining operator $\mathcal{Y}^{\th}$ give nonzero
intertwining operators of types 
$$\fusion{\Ml{\lm}}{\Mt{\pm}}{\Mt{\pm}},\fusion{\Ml{\lm}}{\Mt{\pm}}{\Mt{\mp}}\hbox{ for $\lm\in L^{\circ}$}.$$

\noindent
{\rm(2)} Define $\mathcal{Y}^{\th}\circ{\th}$
by $(\mathcal{Y}^{\th}\circ{\th})(u,z)=\mathcal{Y}^{\th}(\th(u),z)$ for $u\in
V_{L^{\circ}}$.  
Then $\mathcal{Y}^{\th}\circ{\th}$ gives nonzero
intertwining operators for $\M{+}$ of types
$\fusion{\Ml{\lm}}{\Mt{\pm}}{\Mt{}}$.  
Moreover restrictions of $\mathcal{Y}^{\th}$ and
$\mathcal{Y}^{\th}\circ{\th}$ to $\Ml{\lm}\ots\Mt{\pm}$ form a basis of the
vector space $I\fusion{\Ml{\lm}}{\Mt{\pm}}{\Mt{}}$  respectively.
\el
\proof (1) is proved in \cite[Proposition 4.4]{A}.
Clearly $\mathcal{Y}^{\th}\circ{\th}$ gives nonzero intertwining operators of 
types $\fusion{\Ml{\lm}}{\Mt{\pm}}{\Mt{}}$. 
Now we show the second assertion of (2). 
Since $\th\mathcal{Y}^{\th}(u,z)\th(v)=\mathcal{Y}^{\th}(\th(u),z)v$ for 
$u\in\Ml{\lm}$ and $v\in\Mt{\pm}$, 
we have 
$$p_{\pm}((\mathcal{Y}^{\th}\circ{\th})(u,z)v)
=\pm p_{\pm}(\mathcal{Y}^{\th}(u,z)\th(v)),$$
where $p_{\pm}$ is the canonical projection from $\Mt{}$ to $\Mt{\pm}$
respectively.  Hence by Lemma \ref{P1} and (1), we
see that $\mathcal{Y}^{\th}$ and $\mathcal{Y}^{\th}\circ{\th}$ are linearly
independent in the vector spaces $I\fusion{\Ml{\lm}}{\Mt{\pm}}{\Mt{}}$. 
Since the fusion rules of types $\fusion{\Ml{\lm}}{\Mt{\pm}}{\Mt{}}$ are
two by Theorem \ref{T4}, $\mathcal{Y}^{\th}$ and
$\mathcal{Y}^{\th}\circ{\th}$ in fact form a basis of
$I\fusion{\Ml{\lm}}{\Mt{\pm}}{\Mt{}}$.  This proves (2).\qed
\vs

Set $T=T^{1}\ops T^{2}$ the direct sum of the irreducible $\C[L]$-modules
$T^{1}$ and $T^{2}$, and define a linear endomorphism $\psi\in\End T$ by 
$\psi(t_{1})=t_{2},\psi(t_{2})=t_{1}$, where $t_{i}$ is a basis of $T^{i}$ for
$i=1,2$.
For $\lm\in L^{\circ}$, we write $\lm=r\al/2k+m\al$ for $-k+1\leq r\leq k$ and
$m\in\Z$, and define $\psi_{\lm}\in\End T$ by 
$$\psi_{\lm}=e_{m\al}\circ\underbrace{\psi\circ\cdots\circ\psi}_{r}.$$ 
Set $\tilde{\mathcal{Y}}(u,z)=\mathcal{Y}^{\th}(u,z)\ots \psi_{\lm}$ for
$\lm\in L^{\circ}$ and $u\in\Ml{\lm}$, and extend it to $V_{L^{\circ}}$ by
linearity. Then we have following  proposition.

\bp{P4.3}{\rm(1)}  For $\lm\in L^{\circ}$, the linear map $\psi_{\lm}$ has
following properties:
\eqn e_{\beta}\circ\psi_{\lm}=(-1)^{\B\beta,\lm\K}\psi_{\lm}\circ
e_{\beta}=\psi_{\lm+\beta}\hbox{ for all $\beta\in L$}.
\eeqn
\noindent
{\rm(2)} For $a\in\V{}$ and $u\in\Vl{\lm}$,  we have 
\eqn
&& z_{0}^{-1}\delta\left(
{\frac{z_{1}-z_{2}}{z_{0}}}\right)Y^{\th}(a,z_{1})
\tilde{\mathcal{Y}}(u,z_{2})
-\delta\left(
{\frac{z_{2}-z_{1}}{-z_{0}}}\right)\tilde{\mathcal{Y}}(u,z_{2})
Y^{\th}(a,z_{1})\\
&&{ }=\frac{1}{2}\sum_{p=0,1}z_{2}^{-1}\delta\left((-1)^{p}
{\frac{(z_{1}-z_{0})^{1/2}}{z_{2}^{1/2}}}
\right)\tilde{\mathcal{Y}}(Y(\th^{p}(a),z_{0})u,z_{2})
\eeqn
and
$$\frac{d}{dz}\tilde{\mathcal{Y}}(u,z)=\tilde{\mathcal{Y}}(L(-1)u,z).$$
\ep
\proof Since $e_{\al}\circ\psi=-\psi\circ e_{\al}$, we have $e_{m\al}\circ
\psi^{r}=(-1)^{mr}\psi^{r}\circ e_{m\al}$ for $m,r\in\Z$. 
Therefore $\psi_{\lm}$ $(\lm\in L^{\circ})$ satisfies 
$e_{\beta}\circ\psi_{\lm}=(-1)^{\B\beta,\lm\K}\psi_{\lm}\circ e_{\beta}$
and 
$e_{\beta}\circ\psi_{\lm}=\psi_{\lm+\beta}$ for $\beta\in L$. 
This proves (1). Then (2) follows from (\ref{Ja1}), (\ref{Ld1}) and (1).\qed
\vs 

We note that for every quasi-admissible
triple $(\lm,i,j)\in\Qa$, $\psi_{\lm}(T^{i})=T^{j}$. Thus we have 
\bp{P51} Let $(\lm,i,j)\in\Qa$ be an admissible triple. The restriction of 
$\tilde{\mathcal{Y}}$ to $\Vl{\lm}\ots \V{T_{i}}$ gives an intertwining
operator for
$\V{+}$ of type
$\fusion{\Vl{\lm}}{\V{T_{i}}}{\V{T_{j}}}$.
\ep 

Now we have some nonzero intertwining operators by restricting
$\tilde{\mathcal{Y}}$ to irreducible $\V{+}$-modules.
\bp{P32} Fusion rules of following types are nonzero{\rm ;}
\begin{enumerate}
\renewcommand{\labelenumi}{\rm (\roman{enumi})\ \ }
\item
$\fusion{\Vl{\lm_{r}}}{\V{T_{i},\pm}}{\V{T_{j},\pm}},\
\fusion{\Vl{\lm_{r}}}{\V{T_{i},\pm}}{\V{T_{j},\mp}}$ for $r\in\Z$ and
$(\lm_{r},i,j)\in\Qa$,
\item
$\fusion{\V{+}}{\V{T_{i},\pm}}{\V{T_{i},\pm}},\ 
\fusion{\V{-}}{\V{T_{i},\pm}}{\V{T_{i},\mp}}$ for $i\in\{1,2\}$,
\item
$\fusion{\Va{+}}{(\V{T_{1},\pm})'}{\V{T_{1},\pm}},\ 
\fusion{\Va{+}}{(\V{T_{2},\pm})'}{\V{T_{2},\mp}} ,\ 
\fusion{\Va{-}}{(\V{T_{1},\pm})'}{\V{T_{1},\mp}},\ 
\fusion{\Va{-}}{(\V{T_{2},\pm})'}{\V{T_{2},\pm}}$. 
\end{enumerate}
\ep

\proof By Lemma \ref{L3.4} and Proposition \ref{P51}, we see that
$\tilde{\mathcal{Y}}$ gives nonzero intertwining operator of types
$\fusion{\Vl{\lm_{r}}}{\V{T_{i},\pm}}{\V{T_{j},\pm}}$ and 
$\fusion{\Vl{\lm_{r}}}{\V{T_{i},\pm}}{\V{T_{j},\mp}}$ for $r\in\Z$ and
$(\lm_{r},i,j)\in\Qa$. 

Next we shows that fusion rules of types in (ii) and (iii) are nonzero.
By Lemma \ref{L3.4} and Corollary \ref{C8},
$\mathcal{Y}^{\th}(u\pm\th(u),z)v=(\mathcal{Y}^{\th}\pm\mathcal{Y}^{\th}\circ{\th})(u,z)v$
are nonzero for any nonzero $u\in\Ml{\lm}$ ($\lm\in L^{\circ}$) and nonzero
$v\in\Mt{\pm}$. Thus by Proposition \ref{P51}, we see that $\tilde{\mathcal{Y}}$
give nonzero intertwining operators of types 
$$\fusion{\V{+}}{\V{T_{i},\pm}}{\V{T_{i}}},
\fusion{\V{-}}{\V{T_{i},\pm}}{\V{T_{i}}},
\fusion{\Va{+}}{(\V{T_{i},\pm})'}{\V{T_{i}}},
\fusion{\Va{-}}{(\V{T_{i},\pm})'}{\V{T_{i}}}
\hbox{ for $i\in\{1,2\}$}.$$
By the 
definition of $\psi_{\lm}$ $(\lm \in L^{\circ}$), we have
\eqa\psi_{-m\al}=\psi_{m\al},\
\psi_{-(\al/2+m\al)}=e_{-\al}\psi_{\al/2+m\al}\hbox{ for
$m\in\Z$}.\label{60}\eeqa   Since 
$\th\tilde{\mathcal{Y}}(u,z)\th=\mathcal{Y}^{\th}(\th(u),z)\ots \psi_{\lm}$ for $\lm\in L^{\circ}$ and $u\in\Ml{\lm}$, by (\ref{60}) 
we have 
\eqn
&&\th\tilde{\mathcal{Y}}(u,z)\th=\tilde{\mathcal{Y}}(\th(u),z)\hbox{ for
$u\in\V{}$,}\\
&&\th\tilde{\mathcal{Y}}(u,z)\th=e_{\al}\tilde{\mathcal{Y}}(\th(u),z)\hbox{ for
$u\in\Vl{\al/2}$.}
\eeqn
This proves that $\tilde{\mathcal{Y}}$ gives nonzero intertwining operators of
types indicated in (ii) and (iii) of the proposition;  for instance, for $u\in\Va{+}$
and $v\in(\V{T_{2},-})'$, we have 
\eqn
\th\tilde{\mathcal{Y}}(u,z)v&=&e_{\al}\tilde{\mathcal{Y}}(\th(u),z)\th(v)\\
&=&\tilde{\mathcal{Y}}(u,z)v.
\eeqn
Hence $\tilde{\mathcal{Y}}(u,z)v\in\V{T_{2},+}\{z\}$. Thus $\tilde{\mathcal{Y}}$
gives a nonzero intertwining operator of type
$\fusion{\Va{+}}{(\V{T_{2},-})'}{\V{T_{2},+}}$.\qed
\vs

We shall show the following proposition. The proof is given after
Proposition \ref{P17}. 
\bp{P15}$(1)$ For $i,j\in\{1,2\}$, the fusion
rules of types 
$$\fusion{\V{\pm}}{\V{T_{i},\pm}}{\V{T_{j}}}\hbox{ and
}\fusion{\V{\pm}}{\V{T_{i},\mp}}{\V{T_{j}}}$$ 
are zero if $i\neq j$.

\noindent
$(2)$ For $1\leq r \leq k-1$ and $i,j\in\{1,2\}$, the fusion rules of types
$$\fusion{\Vl{\lm_{r}}}{\V{T_{i},\pm}}{\V{T_{j}}}$$
 are zero if
$(-1)^{r+\del_{i,j}+1}\neq1$.

\noindent
$(3)$ For $i,j\in\{1,2\}$, the fusion
rules of types $$\fusion{\Va{\pm}}{\V{T_{i},\pm}}{\V{T_{j}}}
\hbox{ and }\fusion{\Va{\pm}}{\V{T_{i},\mp}}{\V{T_{j}}}$$ 
is zero if $(-1)^{k+\del_{i,j}+1}\neq1$.
\ep

To prove Proposition \ref{P15}, we first show the following proposition.

\bp{P17} Let $W$ be an irreducible $\V{+}$-module and suppose that $W$
contains an $\M{+}$-submodule isomorphic to $\Ml{\lm}$ for some $\lm\in
L^{\circ}$. If
$(\lm,i,j)\in\Pa$ is not a quasi-admissible triple, then fusion rules of types 
$\fusion{W}{\V{T_{i},\pm}}{\V{T_{j}}}$ are zero.
\ep
\proof Let $W$ be an irreducible $\V{+}$-module, 
and suppose that $W$ contains an $\M{+}$-submodule $N$ isomorphic to
$\Ml{\lm}$.  Let $f$ be an $\M{+}$-isomorphism from $\Ml{\lm}$ to $N$. 
Consider an intertwining operator $\mathcal{Y}\in
I_{\V{+}}\fusion{W}{\V{T_{i},\epsilon}}{\V{T_{j}}}$ for
$i,j\in\{1,2\}$ and $\epsilon\in\{\pm\}$. 
We shall prove that $\mathcal{Y}=0$ if $(-1)^{\B\al,\lm\K+\del_{i,j}+1}\neq1.$
  
The restrictions of $\mathcal{Y}$ to $N\ots\V{T_{i},\epsilon}$ gives an intertwining operator for $\M{+}$ of type $\fusion{N}{{\V{T_{i},\epsilon}}}{\V{T_{j}}}.$ 
Set 
\eqn
\overline{\mathcal{Y}}(u,z)
=\phi_{j}^{-1}\mathcal{Y}(f(u),z)\phi_{i}\hbox{ for $u\in\Ml{\lm}$.}
\eeqn
Then $\overline{\mathcal{Y}}$ is an intertwining operator for $\M{+}$ of type
$\fusion{\Ml{\lm}}{\Mt{\epsilon}}{\Mt{}}.$ 
By Lemma \ref{L3.4} (2),  for any $u\in\Ml{\lm}$, 
$\overline{\mathcal{Y}}(u,z)$ is a linear combination of 
$\mathcal{Y}^{\th}(u,z)$ and $\mathcal{Y}^{\th}(\th(u),z)$. By (\ref{Ja1}), we
have 
\eqa
&& z_{0}^{-1}\delta\left(
{\frac{z_{1}-z_{2}}{z_{0}}}\right)\mathcal{Y}^{\th}(E,z_{1})
\mathcal{Y}^{\th}(e_{\pm\lm},z_{2})
-(-1)^{\B\al,\lm\K}z_{0}^{-1}\delta\left(
{\frac{z_{2}-z_{1}}{-z_{0}}}\right)\mathcal{Y}^{\th}(e_{\pm\lm},z_{2})
\mathcal{Y}^{\th}(E,z_{1})\nn\\
&&{}=z_{2}^{-1}\delta\left({\frac{z_{1}-z_{0}}{z_{2}}}
\right)\mathcal{Y}^{\th}(Y(E,z_{0})e_{\pm\lm},z_{2}),\label{Jac1}
\eeqa
where $E=e_{\al}+e_{-\al}\in\V{+}$. Hence one have  
\eqn
(z_{1}-z_{2})^{M}\mathcal{Y}^{\th}(E,z_{1})
\mathcal{Y}^{\th}(e_{\pm\lm},z_{2})
=(-1)^{\B\al,\lm\K}(z_{1}-z_{2})^{M}\mathcal{Y}^{\th}(e_{\pm\lm},z_{2})
\mathcal{Y}^{\th}(E,z_{1})
\eeqn
for a sufficiently large integer $M$, and then 
\eqa
&&(z_{1}-z_{2})^{M}\mathcal{Y}^{\th}(E,z_{1})
\overline{\mathcal{Y}}(e_{\lm},z_{2})
=(-1)^{\B\al,\lm\K}(z_{1}-z_{2})^{M}\overline{\mathcal{Y}}(e_{\lm},z_{2})
\mathcal{Y}^{\th}(E,z_{1}).\label{Ja2}
\eeqa

(\ref{Ja2}) is an identity on $\Mt{\epsilon}$. 
We next derive an identity on $\V{T_{i},\epsilon}$ from (\ref{Ja2}).
Since $e_{\pm\al}\in\C[L]$ act on $\V{T_{i}}$ $(i=1,2)$ as the scalar
$(-1)^{\del_{i,2}}$, we have 
$$
 e_{\pm\al}\phi_{j}\overline{\mathcal{Y}}(u,z)
\phi_{i}^{-1}=(-1)^{\del_{i,j}+1}\phi_{j}
\overline{\mathcal{Y}}(u,z)\phi_{i}^{-1}e_{\pm\al}
$$
for $u\in\Ml{\lm}$. And $Y^{\th}(E,z)$ acts on $\V{T_{i}}$ ($i=1,2$) as
$\mathcal{Y}^{\th}(E,z)\ots e_{\al}$. Hence by (\ref{Ja2}), we have 
\eqa
&&(z_{1}-z_{2})^{M}Y^{\th}(E,z_{1})\mathcal{Y}(f(e_{\lm}),z_{2})\nn\\
&&{\ }=(-1)^{\B\al,\lm\K+\del_{i,j}+1}(z_{1}-z_{2})^{M}\mathcal{Y}(f(e_{\lm}),z_{2})
Y^{\th}(E,z_{1})\label{Ja3}
\eeqa
for a sufficiently large integer $M$. On the
other hand, since $\mathcal{Y}$ is an intertwining operator for $\V{+}$ of
type $\fusion{W}{\V{T_{i},\epsilon}}{\V{T_{j}}},$ Jacobi identity (\ref{JI1}) 
shows that 
\eqa
(z_{1}-z_{2})^{M}Y^{\th}(E,z_{1})\mathcal{Y}(f(e_{\lm}),z_{2})\nn
=(z_{1}-z_{2})^{M}\mathcal{Y}(f(e_{\lm}),z_{2})
Y^{\th}(E,z_{1})\label{Ja4}
\eeqa
for a sufficiently large integer $M$. Therefore by (\ref{Ja3}) and (\ref{Ja4}), if
$(-1)^{\B\al,\lm\K+\del_{i,j}+1}\neq 1$, then 
\eqa(z_{1}-z_{2})^{M}\mathcal{Y}(f(e_{\lm}),z_{2})
Y^{\th}(E,z_{1})u=0\label{Q3}
\eeqa
for a nonzero $u\in\V{T_{i},\epsilon}$ and a sufficiently large integer $M$.
Since there is an integer $n_{0}$ such that $E_{n_{0}}u\neq 0$ and
$E_{n}u=0$ for all $n>n_{0}$, by multiplying $z_{1}^{n_{0}}$ and taking
$\Res_{z_{1}}$  on both side of (\ref{Q3}), we have
$z_{2}^{M}\mathcal{Y}(f(e_{\lm}),z_{2})E_{n_{0}}u=0$. Hence Lemma
\ref{P1} implies that $\mathcal{Y}=0$.\qed
\vs

Now we prove Proposition \ref{P15}.

{\it Proof of Proposition \ref{P15}}.  By the irreducible decompositions
(\ref{dec1})-(\ref{dec3}), we see that
$\Vl{\lm_{r}}$ contains 
$\Ml{\lm_{r}}$ for $1\leq r\leq k-1$, that $\Va{\pm}$ contain an
$\M{+}$-submodule isomorphic to $\Ml{\al/2}$, and that $\V{\pm}$ contain
an $\M{+}$-submodule isomorphic to $\Ml{\al}$. Hence Proposition \ref{P15}
follows from Proposition 
\ref{P17}.\qed

\bp{P21}{\rm(1)} For $i\in\{1,2\}$, fusion rules of types 
$$\fusion{\V{+}}{\V{T_{i},\pm}}{\V{T_{i},\mp}}\hbox{ and }
\fusion{\V{-}}{\V{T_{i},\pm}}{\V{T_{i},\pm}}$$ are zero.

\noindent
{\rm(2)} Fusion rules of types 
$$\fusion{\Va{+}}{(\V{T_{1},\pm})'}{\V{T_{1},\mp}},\ 
\fusion{\Va{+}}{(\V{T_{2},\pm})'}{\V{T_{2},\pm}},\ 
\fusion{\Va{-}}{(\V{T_{1},\pm})'}{\V{T_{1},\pm}},\ 
\fusion{\Va{-}}{(\V{T_{2},\pm})'}{\V{T_{2},\mp}}$$ 
are zero.
\ep
\proof Since $\V{\pm}$ contains the irreducible $\M{+}$-module $\M{\pm}$
respectively and the fusion rules of types  
$\fusion{\M{+}}{\Mt{\pm}}{\Mt{\mp}}\hbox{ and
}\fusion{\M{-}}{\Mt{\pm}}{\Mt{\pm}}$
 are zero by Theorem \ref{T3} (i) and
(ii), (1) follows from Corollary \ref{C8}.

Next we prove that the fusion rules of
types 
$\fusion{\Va{+}}{(\V{T_{1},\pm})'}{\V{T_{1},\mp}}\hbox{ and
}\fusion{\Va{+}}{(\V{T_{2},\pm})'}{\V{T_{2},\pm}}$ are zero. (2) for 
types $\fusion{\Va{-}}{(\V{T_{1},\pm})'}{\V{T_{1},\pm}}\hbox{ and
}\fusion{\Va{-}}{(\V{T_{2},\pm})'}{\V{T_{2},\mp}}$ can be also proved in
the similar way. 

By Proposition \ref{P32} (iii), for $i\in\{1,2\}$ and $\epsilon\in\{\pm\}$,
there exists $\epsilon'\in\{\pm\}$ such that the fusion rule of type
$\fusion{\Va{+}}{(\V{T_{i},\epsilon})'}{\V{T_{i},\epsilon'}}$ is nonzero.  Let
$\{\tau,\epsilon'\}=\{\pm\}$. Then we have to prove that the fusion rule of type
$\fusion{\Va{+}}{(\V{T_{i},\epsilon})'}{\V{T_{i},\tau}}$ is zero.  To show this, we
prove that the canonical projection
$$\fusion{\Va{+}}{(\V{T_{i},\epsilon})'}{\V{T_{i}}}\to\fusion{\Va{+}}{(\V{T_{i},
\epsilon})'}{\V{T_{i},\epsilon'}},\
\mathcal{Y}\mapsto p_{\epsilon'}\circ\mathcal{Y}$$ is injective, where $p_{\pm}$
are the canonical projection from $\V{T_{i}}$ to $\V{T_{i},\pm}$ respectively  and
$p_{\pm}\circ\mathcal{Y}$ are intertwining operators defined by
$(p_{\pm}\circ\mathcal{Y})(u,z)v=p_{\pm}(\mathcal{Y}(u,z)v)$ for $u\in\Va{+}$
and $v\in(\V{T_{i},\epsilon})'$. To prove this, it is enough to prove that arbitrary
nonzero intertwining operator $\mathcal{Y}$ of type
$\fusion{\Va{+}}{(\V{T_{i},\epsilon})'}{\V{T_{i}}}$ satisfies
\eqa
\th\mathcal{Y}(e_{\al/2}+e_{-\al/2},z)\th
=(-1)^{\del_{i,2}}\mathcal{Y}(e_{\al/2}+e_{-\al/2},z).\label{RT1}
\eeqa
Actually if $\mathcal{Y}$ is a nonzero intertwining operator of the indicated type which satisfies (\ref{RT1}), then $p_{\epsilon'}(\mathcal{Y}(e_{\al/2}+e_{-\al/2},z)v)=\mathcal{Y}(e_{\al/2}+e_{-\al/2},z)v$ for $v\in(\V{T_{i},\epsilon})'$ and then $p_{\epsilon'}\circ\mathcal{Y}$ is nonzero by Corollary \ref{C8}.

Let $\V{T_{j}}=(\V{T_{i}})'$, and let $\Va{+}[0]$ be as
of (\ref{VA1}). Then $\mathcal{Y}$ gives an intertwining operator for
$\M{+}$ of type $\fusion{\Va{+}[0]}{\V{T_{j},\epsilon}}{\V{T_{i}}}.$ Thus by
Lemma \ref{L3.4} (2), we see that 
$I_{\M{+}}\fusion{\Va{+}[0]}{\V{T_{i},\epsilon}}{\V{T_{j}}}$ is spanned by
intertwining operators $\mathcal{Y}^{\pm}$ defined by
\eqa
\mathcal{Y}^{\pm}(u,z)=\phi_{i}\mathcal{Y}^{\th}(\phi_{\pm\al/2}(u),z)
\phi_{j}^{-1}\hbox{ for $u\in \Va{+}[0]$}.\label{TW1}
\eeqa
Hence there exist constants $c_{1},c_{2}\in\C$ such that 
\eqa
\mathcal{Y}(u,z)=c_{1}\mathcal{Y}^{+}(u,z)+c_{2}\mathcal{Y}^{-}(u,z)
\label{TW2}
\eeqa
for all $u\in\Va{+}[0]$. Now for $\beta\in\h$, set 
\eqn
\exp\left(\sum_{n=0}^{\infty}\frac{\beta(-n)}{n}
z^{n}\right)=\sum_{n=0}^{\infty}p_{n}(\beta)z^{n}\in(\End V_{L^{\circ}})[[z]].
\eeqn
Then we have $E_{0}(e_{\frac{\al}{2}}+e_{-\frac{\al}{2}})=p_{k-1}(\al)
e_{\frac{\al}{2}}+p_{k-1}(-\al)e_{-\frac{\al}{2}}\in \Va{+}[0]$, and hence 
\eqn
\phi_{\frac{\al}{2}}(E_{0}(e_{\frac{\al}{2}}+e_{-\frac{\al}{2}}))
=p_{k-1}(\al)e_{\frac{\al}{2}},\ \ 
\phi_{-\frac{\al}{2}}(E_{0}(e_{\frac{\al}{2}}+e_{-\frac{\al}{2}}))
=p_{k-1}(-\al)e_{-\frac{\al}{2}}.
\eeqn
Thus by (\ref{TW1}) and (\ref{TW2}), we have
\eqa
&&[E_{0},\mathcal{Y}(e_{\frac{\al}{2}}+
e_{-\frac{\al}{2}},z)]\nn\\
&&\ \ = \mathcal{Y}(E_{0}(e_{\frac{\al}{2}}+e_{-\frac{\al}{2}}),z)\nn\\
&&\ \ =\phi_{i}(c_{1}\mathcal{Y}^{\th}(p_{k-1}(\al)e_{\frac{\al}{2}},z)
+c_{2}\mathcal{Y}^{\th}(p_{k-1}(-\al)e_{-\frac{\al}{2}},z))
\phi_{j}^{-1}.\label{TW3}
\eeqa
On the other hand, (\ref{Jac1}) shows that
\eqn 
[E_{0}, \phi_{i}\mathcal{Y}^{\th}(e_{\pm\frac{\al}{2}},z)\phi_{j}^{-1}]
&=&e_{\al}\phi_{i}\mathcal{Y}^{\th}(E_{0}(e_{\pm\frac{\al}{2}}),z)\phi_{j}^{-1}\\
&=&(-1)^{\del_{i,2}}\phi_{i}\mathcal{Y}^{\th}(p_{k-1}(\mp\al)
e_{\mp\frac{\al}{2}},z)\phi_{j}^{-1}.
\eeqn
Hence by (\ref{TW1}) and (\ref{TW2}) again, we have
\eqa
&&[E_{0},\mathcal{Y}(e_{\frac{\al}{2}}
+e_{-\frac{\al}{2}},z)]\nn\\
&&\ \ =
(-1)^{\del_{i,2}}c_{1}([E_{0},\mathcal{Y}^{+}(E_{0}(e_{\frac{\al}{2}}
+e_{-\frac{\al}{2}}),z)]
+c_{2}[E_{0},\mathcal{Y}^{-}(E_{0}(e_{\frac{\al}{2}}+e_{-\frac{\al}{2}}),z)])
\nn\\
&&\ \ =
(-1)^{\del_{i,2}}\phi_{i}(c_{1}\mathcal{Y}^{\th}(p_{k-1}(-\al)
e_{-\frac{\al}{2}},z)
+c_{2}\mathcal{Y}^{\th}(p_{k-1}(\al)e_{\frac{\al}{2}},z))
\phi_{j}^{-1}.\label{TW4}
\eeqa
Subtracting (\ref{TW3}) from (\ref{TW4}) gives the identity
\eqn
(c_{1}-(-1)^{\del_{i,2}}c_{2})\phi_{i}(\mathcal{Y}^{\th}(p_{k-1}(\al)
e_{\frac{\al}{2}},z)
-(-1)^{\del_{i,2}}\mathcal{Y}^{\th}(p_{k-1}(-\al)e_{-\frac{\al}{2}},z))
\phi_{j}^{-1}=0.
\eeqn
 Then Lemma \ref{L3.4} shows that $c_{1}=(-1)^{\del_{i,2}}c_{2}$. Since
$\th\mathcal{Y}^{\pm}(u, z)\th=\mathcal{Y}^{\mp}(u,z)$ for $u\in\Va{+}[0]$, 
we have (\ref{RT1}).\qed
\vs

Now the following proposition follows from Proposition \ref{P9}, Proposition 
\ref{P49}, Proposition \ref{P32}, Proposition \ref{P15} and Proposition \ref{P21}.

\bp{MP2} Let $W^{1},W^{2}$ and $W^{3}$ be irreducible
$\V{+}$-modules and suppose that one of them is of twisted type. Then the
fusion rule $N_{W^{1}W^{2}}^{W^{3}}$ is zero or one.
Assume that $W^{1}$ is twisted type module, then the fusion rule
$N_{W^{1}W^{2}}^{W^{3}}$ is one if and only if $W^{i}$ $(i=1,2,3)$ satisfy the
following cases{\rm :} 

\noindent
{\rm(i)} $W^{1}=(\V{T_{1},+})'$ and the pair $(W^{2},W^{3})$ is one of pairs
\eqn &&(\V{\pm},(\V{T_{1},\pm})'),\ (\V{T_{1},\pm},\V{\pm}),
\ (\Va{\pm},\V{T_{1},\pm}),\ ((\V{T_{1},\pm})',
(\Va{\pm})'),\\ 
&&(\Vl{\lm_{r}},
(\V{T_{1},\pm})')\hbox{ and }(\V{T_{1},\pm},\Vl{\lm_{r}})\hbox{ for $1\leq
r\leq k-1$ and $r$ is even},\\ 
&&(\Vl{\lm_{r}},(\V{T_{2},\pm})')\hbox{ and
}(\V{T_{2},\pm},\Vl{\lm_{r}})\hbox{ for $1\leq
r\leq k-1$ and $r$ is odd}.
\eeqn

\noindent
{\rm(ii)} $W^{1}=(\V{T_{1},-})'$ and the pair $(W^{2},W^{3})$ is one of
pairs
\eqn &&(\V{\pm},(\V{T_{1},\mp})'),\ (\V{T_{1},\pm},\V{\mp}),\ 
(\Va{\pm},\V{T_{1},\mp}),\ ((\V{T_{1},\pm})',(\Va{\mp})'),\\ &&(\Vl{\lm_{r}},
(\V{T_{1},\pm})')\hbox{ and }(\V{T_{1},\pm},\Vl{\lm_{r}})\hbox{ for $1\leq
r\leq k-1$ and $r$ is even},\\ 
&&(\Vl{\lm_{r}},(\V{T_{2},\pm})')\hbox{ and
}(\V{T_{2},\pm},\Vl{\lm_{r}})\hbox{ for $1\leq
r\leq k-1$ and $r$ is odd}.
\eeqn

\noindent
{\rm(iii)} $W^{1}=(\V{T_{2},+})'$ and the pair $(W^{2},W^{3})$ is one of
pairs
\eqn &&(\V{\pm},(\V{T_{2},\pm})'),\ (\V{T_{2},\pm},\V{\pm}),\ 
(\Va{\pm},\V{T_{2},\mp}),\ ((\V{T_{2},\pm})',(\Va{\mp})'),\\ 
&&(\Vl{\lm_{r}},(\V{T_{2},\pm})')\hbox{ and }(\V{T_{2},\pm},\Vl{\lm_{r}})\hbox{
for $1\leq r\leq k-1$ and $r$ is even},\\ 
&&(\Vl{\lm_{r}},(\V{T_{1},\pm})')\hbox{ and
}(\V{T_{1},\pm},\Vl{\lm_{r}})\hbox{ for $1\leq
r\leq k-1$ and $r$ is odd}.
\eeqn

\noindent
{\rm(iv)} $W^{1}=(\V{T_{2},-})'$ and the pair $(W^{2},W^{3})$ is one of
pairs
\eqn &&(\V{\pm},(\V{T_{2},\mp})'),\ (\V{T_{2},\pm},\V{\mp}),\ 
(\Va{\pm},\V{T_{2},\pm}),\ ((\V{T_{2},\pm})',(\Va{\pm})'),\\
&&(\Vl{\lm_{r}},(\V{T_{2},\pm})')\hbox{ and }(\V{T_{2},\pm},\Vl{\lm_{r}})\hbox{
for $1\leq r\leq k-1$ and $r$ is even},\\  &&(\Vl{\lm_{r}},(\V{T_{1},\pm})')\hbox{
and }(\V{T_{1},\pm},\Vl{\lm_{r}})\hbox{ for $1\leq r\leq k-1$ and $r$ is odd}.
\eeqn
\ep

Consequently Theorem \ref{T4} follows from Proposition \ref{P9}, Proposition \ref{MP1} and Proposition
\ref{MP2}.

\begin{center} {\bf ACKNOWLEDGMENT}\end{center}

I would like to thank Professor Kiyokazu Nagatomo for helpful advice.


\end{document}